\documentclass[10pt]{amsart}
\usepackage[letterpaper,textwidth=6.15in,textheight=8.65in,centering]{geometry}

\usepackage[T1]{fontenc}
\usepackage{microtype}
\usepackage{amsmath,amssymb,mathtools}

\usepackage{enumitem}
\usepackage{hyperref}
\hypersetup{
  pdftitle={Stability conditions on twisted derived categories},
  pdfauthor={Zhirui Li},
  colorlinks=false,
  citebordercolor={0 1 0},
  linkbordercolor={1 0 0},
  urlbordercolor={0 1 1},
  pdfborder={0 0 1}
}
\numberwithin{equation}{section}
\usepackage{xcolor}

\newtheorem{theorem}{Theorem}[section]
\newtheorem{proposition}[theorem]{Proposition}
\newtheorem{lemma}[theorem]{Lemma}
\newtheorem{corollary}[theorem]{Corollary}
\theoremstyle{definition}
\newtheorem{definition}[theorem]{Definition}

\theoremstyle{remark}
\newtheorem{remark}[theorem]{Remark}

\newcommand{\C}{\mathbb C}
\newcommand{\Q}{\mathbb Q}
\newcommand{\R}{\mathbb R}
\newcommand{\Z}{\mathbb Z}
\newcommand{\PP}{\mathbb P}
\newcommand{\OO}{\mathcal O}
\newcommand{\CC}{\mathcal C}
\newcommand{\DD}{\mathcal D}
\newcommand{\HH}{\mathcal H}
\newcommand{\cA}{\mathcal A}
\newcommand{\BB}{\mathcal B}
\newcommand{\cP}{\mathcal P}
\newcommand{\cQ}{\mathcal Q}
\newcommand{\cR}{\mathcal R}
\newcommand{\Db}{\mathrm D^b}
\newcommand{\Dqc}{\mathrm D_{\mathrm{qc}}}
\newcommand{\Perf}{\mathrm{Perf}}
\newcommand{\Coh}{\operatorname{Coh}}
\newcommand{\Ind}{\operatorname{Ind}}
\newcommand{\Hom}{\operatorname{Hom}}

\newcommand{\Cone}{\operatorname{Cone}}
\newcommand{\im}{\operatorname{im}}
\newcommand{\rk}{\operatorname{rk}}
\newcommand{\ch}{\operatorname{ch}}
\newcommand{\trd}{\operatorname{trd}}
\newcommand{\Knum}{K_{\mathrm{num}}}
\newcommand{\Stab}{\operatorname{Stab}}
\newcommand{\id}{\mathrm{id}}
\newcommand{\eps}{\varepsilon}
\newcommand{\phip}{\phi^{+}}
\newcommand{\phim}{\phi^{-}}

\title[Stability conditions on twisted derived categories]
{Stability conditions on twisted derived categories I: Existence and full numerical support}

\author{Zhirui Li}
\address{Qiuzhen College, Tsinghua University,
Beijing 100084, China}
\email{lizhirui122@gmail.com}
\date{November 2026}

\begin{document}
\begin{abstract}
We prove that for every nonempty smooth projective variety $X$ over an arbitrary field and every Brauer class $\alpha$, 
the twisted derived category $\Db(X,\alpha)$ admits a Bridgeland stability condition with full numerical support.
Based on Chunyi Li's construction and its extension to arbitrary fields,
we establish uniform phase control under tensoring by an Azumaya algebra and calibrate the resulting slicing to the induced central charge.
We further show that the constructed stability condition is geometric and satisfies a mass--Hom bound, 
which persists throughout its connected component.
\end{abstract}
\maketitle
\setcounter{tocdepth}{1}
\tableofcontents

\section{Introduction}
Motivated by Douglas's notion of $\Pi$-stability for D-branes in string theory \cite{Douglas01,Douglas02},
Bridgeland introduced stability conditions on triangulated categories in \cite{Bri}.
His fundamental deformation theorem \cite[Theorem~1.2]{Bri} established the local description of spaces of stability conditions in terms of their central charges.
In particular,
numerical stability spaces with finite-rank numerical Grothendieck groups are finite-dimensional complex manifolds by \cite[Corollary~1.3]{Bri}.

Constructing Bridgeland stability conditions on \(\Db(X)\) for a scheme $X$ thus became a central problem in algebraic geometry. 
The case of smooth projective curves was understood through slope stability by \cite{Bri,Macri07}.
Bridgeland constructed stability conditions on K3 surfaces and described a distinguished connected component of their stability manifolds in \cite{BriK3}.
The construction via tilted hearts was subsequently developed for arbitrary smooth projective complex surfaces by Arcara--Bertram \cite{AB13}, 
using classical Bogomolov--Gieseker inequality and the Hodge index theorem.

For threefolds,
Bayer,
Macr\i` and Toda proposed a double-tilt construction \cite{BMT14},
relating the existence of Bridgeland stability conditions to Bogomolov--Gieseker-type inequalities involving $\ch_3$ of tilt-semistable objects.

Subsequent work established stability conditions on several classes of threefolds,
including projective space \cite{Macri14},
smooth quadric threefolds \cite{Schmidt14},
Fano threefolds of Picard rank one \cite{LiFano19}, 
abelian threefolds and certain Calabi--Yau threefolds of abelian or Kummer type \cite{BMS16},
and smooth quintic threefolds \cite{LiQuintic19}.
Further examples include Calabi--Yau double and triple solids \cite{Koseki22} and broader classes of Calabi--Yau threefolds \cite{FKLR25}.

Beyond existence,
an important property is \emph{geometricity}, 
which requires all skyscraper sheaves of points to be stable of the same phase \cite{BriK3,LiQuintic19,FKLR25}.
Further questions concern the full numerical support property \cite{LLLMPSZ,Cheng} and the geometry of moduli spaces of semistable objects \cite{BM14,PT19,BLMNPS}.

In January 2026,
Li \cite{Li} proved that for any general smooth projective variety $X$ over $\C$,
$\Stab(\Db(X))\neq\varnothing$.
Based on stability conditions on products of elliptic curves by \cite{Liu21},
Li constructed stability conditions on $\PP^n$ via finite quotients of $E^n$ and then induced them on $X\hookrightarrow\PP^n$ by restriction.
This approach does not rely on the higher-dimensional Bogomolov--Gieseker type inequalities used by earlier constructions.

In July 2026,
Li et al.\ \cite{LLLMPSZ} extended the construction of stability conditions in \cite{Li} to projective schemes over arbitrary fields and to projective families, 
based on the relative theory of \cite{BLMNPS}.
They constructed geometric stability conditions, 
established mass--Hom bounds throughout a distinguished component \cite[Theorem~1.3]{LLLMPSZ},
and obtained proper good moduli spaces in characteristic zero \cite[Theorem~1.1]{LLLMPSZ} via \cite[Theorem~21.24]{BLMNPS}.

In \cite[Section~1.7]{LLLMPSZ},
the authors raised three directions for further research: 
\begin{itemize}[nosep]
    \item Projectivity of the constructed proper good moduli spaces;
    \item Full numerical support on $\Knum(X)$;
    \item Stability conditions on other categories, including $\Db(X,\alpha)$.
\end{itemize}
The second question was subsequently solved by Cheng \cite{Cheng} in August 2026.
The third concerns smooth proper Deligne--Mumford stacks, 
twisted derived categories and semiorthogonal components of derived categories, 
with important examples of the latter constructed in \cite{BLMS23,PPZ22}.

Twisted derived categories $\Db(X,\alpha)$ arise naturally from moduli problems for sheaves \cite{Caldararu00, Yoshioka06, Lieblich07} and play an important role in Fourier--Mukai equivalences \cite{Caldararu00,HS06}.
Earlier constructions include generic twisted K3 surfaces \cite{HMS08}, 
arbitrary smooth projective surfaces \cite[Theorem~12.2]{HP24} and twisted abelian threefolds \cite[Theorems~12.6 and~12.8]{HP24} in characteristic zero.
However,
these results do not provide a general existence theorem with full numerical support in arbitrary dimension and characteristic.

We address this problem for all smooth projective varieties over arbitrary fields and all Brauer classes,
as stated in the following theorem.

\begin{theorem}\label{thm:main}
Let $X$ be a nonempty smooth projective scheme over an arbitrary field $k$,
and let $\alpha\in\operatorname{Br}(X)$.
Then $\Db(X,\alpha)$ admits a Bridgeland stability condition satisfying the support property with respect to its full numerical Grothendieck group $\Knum(\Db(X,\alpha))$.
\end{theorem}

The constructed stability condition is also geometric and satisfies a mass--Hom bound, 
and the mass--Hom bound persists throughout its connected component; 
see Proposition~\ref{prop:properties-mass-Hom} and Corollary~\ref{cor:properties-component-mass-hom},
together with the geometricity result in Section~\ref{subsec:geometricity}.

Although $\Db(X,\alpha)$ can be realized as an admissible subcategory of the derived category of a Brauer--Severi variety \cite{BerghSchnuerer},
a stability condition on the ambient category need not restrict to $\Db(X,\alpha)$, 
since the Harder--Narasimhan factors of an object in the subcategory need not remain in it.

We briefly outline the proof of Theorem~\ref{thm:main}.
Choose an Azumaya algebra $\cA$ representing $\alpha$,
and write 
$$\CC=\Db(X),\qquad \DD=\Db(X,\cA).$$
Let $L:\CC\rightleftarrows\DD:U$ be the free--forgetful adjunction.
Starting from a suitable rational stability condition $\sigma=(Z,\cP)$ constructed in \cite[Theorem~6.2]{LLLMPSZ},
we establish uniform phase control under tensoring by $\cA$ in Proposition~\ref{lem:bundle}:
$$
UL(\cP(\phi)) \subseteq \cP([\phi-\eps,\phi+\eps]).
$$
This estimate is recorded in \eqref{eq:fixed-error}. 
While $W=Z\circ U$ is the natural candidate for a central charge on $\DD$, the main technical challenge is to construct a locally finite slicing $\cR$ such that $(W,\cR)$ is the stability we want.
Here the key technique is construct a locally finite slicing $\cQ$ whose phase are uniformly close to the phase determined by $W$ first, see Lemma~\ref{lem:large-bands} to Corollary~\ref{cor:induced-charge}; so that $\cQ$ can be subsequently modified to a slicing $\cR$ compatible to $W$ in Theorem~\ref{thm:calibration}.

For each $c\in\R$,
we generate a $t$-structure $(\widehat{\cA}_c,\widehat{\BB}_c)$ on $\Dqc(X,\cA)$ from $L(\cP(>c))$,
whose coaisle is $U^{-1}\widehat{\cP}(\leq c)$.
The key step is Proposition~\ref{prop:coherence},
which proves that $(\widehat{\cA}_c,\widehat{\BB}_c)$ restricts to a bounded $t$-structure $(\cA_c,\BB_c)$ on $\DD$.
The proof uses noetherian hearts of $\sigma$ to ensure that their truncations preserve bounded coherent complexes.

For this family of $t$-structures,
we define $\cQ(\phi)$ by \eqref{eq:twisted-slicing}.
Proposition~\ref{prop:finite-transitions} shows that $\cQ$ is a locally finite slicing on $\DD$,
satisfying
$$U(\cQ(\phi))\subseteq\cP([\phi-\eps,\phi]).$$
Adapting the local phase-calibration argument of \cite[Section~7]{Bri},
Theorem~\ref{thm:calibration} produces a locally finite slicing $\cR$ compatible with $W=Z\circ U$,
and yields a locally finite pre-stability condition $(W,\cR)$.

Finally,
Proposition~\ref{lem:K0-comparison} shows that the forgetful functor induces an injection 
$$
U_{\mathrm{num}}:\Knum(\DD)\hookrightarrow \Knum(\CC),
$$
which becomes an isomorphism after tensoring with $\Q$.
Combining the mass estimate in Lemma~\ref{lem:mass estimate} with \cite[Corollary~2.3]{Cheng},
we obtain the full support property for $\tau=(W,\cR)$.
The disconnected case follows by treating the connected components separately and applying Morita equivalence.

In Section~\ref{sec:properties},
we establish two further properties of $\tau$.
First,
using the compatibility of our construction with scalar extension stated in Proposition~\ref{lem:properties-base-change},
we prove that $\tau$ is geometric in Proposition~\ref{prop:properties-geometric}.
Second,
we establish the mass--Hom bound for $\tau$ in Proposition~\ref{prop:properties-mass-Hom} and show that it persists throughout its connected component in Corollary~\ref{cor:properties-component-mass-hom}.

\begin{remark}\label{rem:good-moduli}
Assume $\operatorname{char}(k)=0$;
then the category $\DD=\Perf(X,\cA)$ is smooth,
proper,
and idempotent complete.
By Proposition~\ref{prop:properties-mass-Hom} and Corollary~\ref{cor:properties-component-mass-hom},
every stability condition $\rho$ in the connected component $\Sigma$ containing $\tau$ satisfies a mass--Hom bound.

Let $\rho=(V,\mathcal{S})\in\Sigma$ and $N=\Knum(\DD)$.
By \cite[Theorem~12.17 and Proposition and Definition~12.15]{BLMNPS},
for every field extension $K/k$,
the canonical extension $\rho_K=(V_K,\mathcal{S}_K)$ has numerical class map
$$
v_K:\Knum(\DD_K)\to \Lambda_k,\qquad v_K(E)=\eta^{\vee}_{K/k}(E),
$$
where $\Lambda_k=\operatorname{Im}(\eta^{\vee}_{\bar k/k})\subset N\otimes_{\Z}\Q$.
The identity $\chi_k(T,v_K(E))=\chi_K(T_K,E)$ shows that these classes are universally locally constant in perfect families,
by proper base change.
Together with the uniform support property under scalar extension,
this gives a locally constant stability condition in the sense of \cite[Example~2.29]{HLR}.

Let $H=LG$ be the classical generator of $\DD$ proved in Lemma~\ref{lem:azumaya-compacts}.
For every finite extension $K/k$ and $E\in\DD_K$,
the adjunction formula and restriction of scalars give
\[
\begin{aligned}
[K:k]\dim_K\Hom_{\DD_K}(H_K,E)
&=\dim_k\Hom_{\DD}(H,p_*E)\\
&\leq C_{H,\rho}\,m_\rho(p_*E)\\
&=C_{H,\rho}[K:k]m_{\rho_K}(E).
\end{aligned}
\]
The constant is uniform over all field extensions,
which verifies condition (1) of \cite[Theorem~2.35]{HLR}. 
Moreover,
the generator boundedness condition holds by \cite[Remark~2.36]{HLR}.
Consequently,
\cite[Theorem~2.35]{HLR} implies that for every $\mathbf v\in\Lambda_k$ and $a\in\R$, 
the locus of $\rho_K$-semistable objects of class $\mathbf v$ in $\mathcal S_K((a,a+1])$ defines an open derived substack $\mathfrak M_{\rho}^{\mathrm{ss}} (\mathbf v;a)$ with a canonical morphism
\[
\pi: \mathfrak M_{\rho}^{\mathrm{ss}}
(\mathbf v;a)_{\mathrm{cl}}
\longrightarrow
M_{\rho}^{\mathrm{ss}}(\mathbf v;a),
\]
where $M_{\rho}^{\mathrm{ss}}(\mathbf v;a)$ is a proper good moduli space over $k$.

For the untwisted projective families in \cite[Theorem~1.1]{LLLMPSZ},
the existence of proper good moduli spaces in characteristic zero follows from \cite[Theorem~21.24]{BLMNPS}.
Our application of \cite[Theorem~2.35]{HLR} gives an analogous result for twisted derived categories over a field, 
without establishing a relative theory for twisted projective families.

In both settings,
properness concerns the good moduli spaces rather than the semistable stacks, 
and projectivity is not established in general. 
In particular,
our result does not solve the first question raised in \cite[Section~1.7]{LLLMPSZ}.
\end{remark}

\medskip
\noindent\textbf{Statements.}

The author has no relevant financial or non-financial interests to disclose.
Data sharing is not applicable to this article as no datasets were used.

\medskip
\noindent\textbf{AI disclosure.}

The initial construction of stability conditions on twisted derived categories was first proposed by ChatGPT 6 Astra, 
an AI system developed by OpenAI. 
After independently verifying this construction and its proofs, 
the author further asked the same AI system whether the resulting stability conditions satisfied additional properties, such as geometricity and mass--Hom bound.  
In response, the AI provided mathematical arguments for these 2 properties, and the moduli-space application presented in Remark~\ref{rem:good-moduli}.

ChatGPT 6 Astra and Gemini 4 Argon were also used to assist with checking mathematical arguments and improving the exposition. 
The author has independently examined and carefully verified all mathematical statements and proofs in this manuscript, and has rewritten almost all of the proofs originally generated by AI.
The author takes full responsibility for the mathematical content, references and any remaining errors in the final manuscript.

\medskip
\noindent\textbf{Acknowledgement.}

The author is grateful to Chunyi Li for several discussions on related topics in January and April 2026. 
These discussions did not involve the specific problem in this paper or the construction method developed here. 
The author also thanks Ziqi Liu for his valuable comments and suggestions on an earlier version of the manuscript, 
particularly for suggesting a clearer explanation of the motivation behind the construction developed 
from Lemma~\ref{lem:large-bands} through Corollary~\ref{cor:induced-charge}. 
The author would also like to thank his supervisor, Zhiwei Zheng, for his continued support throughout this project.

After completing the first draft of this manuscript on November 10, 2026, the author knew from Chunyi Li and Ziqi Liu that at least two other people had independently obtained similar results with the assistance of AI tools, 
including Ziqi Liu in September 2026. The results and arguments presented in this paper were developed independently, 
and the author was unaware when preparing the first manuscript.

\section{Preliminaries}\label{sec:prelim}

\subsection{Stability conditions and scalar extension}
\label{subsec:prelim-stability}

We recall the conventions for Bridgeland stability conditions and the basic results on numerical support and scalar extension used throughout the paper.

Let $\mathcal T$ be an essentially small triangulated category. 
For full subcategories $\mathcal S,\mathcal T'\subseteq\mathcal T$,
write $\mathcal S*\mathcal T'$ for their extension product, 
consisting of objects $E$ fitting into a distinguished triangle
\[
S\longrightarrow E\longrightarrow T'\longrightarrow S[1],
\qquad S\in\mathcal S,\quad T'\in\mathcal T'.
\]
The operation $*$ is associative by the octahedral axiom.

An abelian category is called \textit{noetherian} (respectively,
\textit{artinian}) if every ascending (respectively,
descending) chain of subobjects stabilizes. 
It is of \textit{finite length} if it is both noetherian and artinian. 
For the quasi-abelian categories considered below,
finite length is understood with respect to strict subobjects.

We now recall the notion of a slicing.

\begin{definition}\label{prelim:slicing}
A \emph{slicing} $\cP$ on $\mathcal{T}$ consists of full additive subcategories $\cP(\phi)\subset\mathcal T$,
indexed by $\phi\in\R$,
satisfying
$$
\cP(\phi+1)=\cP(\phi)[1],\qquad\Hom(\cP(\phi_1),\cP(\phi_2))=0\quad\text{for every $\phi_1>\phi_2$,}
$$
and such that every nonzero object $E$ admits a finite \emph{Harder--Narasimhan filtration} with nonzero factors $F_j\in \cP(\phi_j)$ of strictly decreasing phases $\phi_1>\cdots>\phi_n$.

We write $$
\phi^+_{\cP}(E)=\phi_1,\qquad \phi_{\cP}^{-}(E)=\phi_n.
$$
For an interval $I\subseteq\R$,
the category $\cP(I)$ is the finite extension closure of the slices with phases in $I$,
including the zero object.
The nonzero objects of $\cP(\phi)$ are called $\cP$-\textit{semistable} objects of phase $\phi$, and its simple objects are called \textit{stable}.

A slicing is \emph{locally finite} if there exists $\eta>0$ such that $\cP((\phi-\eta,\phi+\eta))$ is a finite length quasi-abelian category,
where quasi-abelian is defined in \cite[Definition~4.1]{Bri}.
\end{definition}
We next recall the numerical invariants used to define central charges and the support property for stability conditions.

The Grothendieck group $K_0(\mathcal T)$ is generated by classes of objects,
satisfying $[E]=[F]+[G]$ for every distinguished triangle $F\to E\to G\to F[1]$.
If $\mathcal T$ is \textit{proper} over a field $k$,
its Euler pairing is defined by
\[
\chi(E,F)=\sum_{j\in\Z}(-1)^j\dim_k\Hom_{\mathcal T}(E,F[j]).
\]
Here proper means for every $E,F\in\mathcal{T}$,
$\sum_{j\in\Z}\dim_k\Hom_{\mathcal{T}}(E,F[j])<\infty$.
The numerical Grothendieck group is $\Knum(\mathcal T)=K_0(\mathcal T)/\operatorname{rad}_R\chi$,
where
\[
\operatorname{rad}_R\chi=
\left\{v\in K_0(\mathcal T):\chi(u,v)=0 \text{ for every }u\in K_0(\mathcal T)
\right\}.
\]
We now recall Bridgeland stability conditions and the support property.
\begin{definition}
Let $\Lambda$ be a finite-rank free abelian group,
and fix a homomorphism $v:K_0(\mathcal{T})\to\Lambda$.
A \emph{pre-stability condition} is a pair $\sigma=(Z,\cP)$,
where $\cP$ is a slicing and $Z:\Lambda\to\C$ is a homomorphism satisfying
$$
Z(v(E))\in\R_{>0}e^{i\pi\phi}\qquad \text{for every $0\neq E\in\cP(\phi)$.}
$$
We say that $\sigma$ satisfies the \textit{support property} if there exists a norm $\|\cdot\|$ on $\Lambda_{\R}=\Lambda\otimes_{\Z}\R$ and a constant $C_{\sigma}>0$ such that
\begin{equation}\label{eq:support}
\|v(E)\|\leq C_{\sigma} |Z(v(E))|\qquad \text{for every nonzero $\sigma$-semistable object $E$.}    
\end{equation}
A \textit{Bridgeland stability condition} with respect to $\Lambda$ is a pre-stability condition satisfying the support property.
We denote the space of such conditions by $\Stab_{(\Lambda,v)}(\mathcal{T})$.

When $\mathcal T$ is proper and $\Lambda=\Knum(\mathcal T)$,
with $v$ being the natural quotient map $K_0(\mathcal{T})\to\Knum(\mathcal{T})$, we call this the \emph{full numerical support property}.
By \cite[Theorem~1.2]{Bri} and \cite[Theorem~12.11]{BLMNPS},
$\Stab_{(\Lambda,v)}(\mathcal T)$ is a complex manifold,
and the central-charge map $(Z,\cP)\to Z$ is a local homeomorphism from $\Stab_{(\Lambda,v)}(\mathcal{T})$ to $\Hom_{\Z}(\Lambda,\C)$.
\end{definition}

\paragraph{\textbf{\textit{Cuts, hearts, and exact sequences.}}} 
A slicing determines a family of bounded \(t\)-structures.
We recall the corresponding conventions for cuts and hearts.
\begin{definition}\label{def:prelim-cuts-hearts}
A $t$-structure on $\mathcal T$ is described by an orthogonal pair $(\cA,\BB)$ satisfying
\[
\cA[1]\subseteq\cA,
\qquad
\BB=\cA^\perp,
\]
such that every $E\in\mathcal T$ admits a truncation triangle $E^+\to E\to E^- \stackrel{+1}{\to}$ with $E^+\in\cA$ and $E^-\in\BB$.
In the usual indexing convention,
\[
\mathcal T^{\leq0}=\cA,
\qquad
\mathcal T^{\geq0}=\BB[1],
\]
and the heart is $\cA\cap\BB[1]$.

For a slicing $\cP$ and $d\in\R$,
the associated cut is $(\cP(>d),\cP(\leq d))$,
with bounded heart
\[
\HH_d=\cP((d,d+1]),
\qquad
\mathcal T_d^{\leq0}=\cP(>d),
\qquad
\mathcal T_d^{\geq0}=\cP(\leq d+1).
\]
\end{definition}

We use the convention $\Coh(X)[0,+\infty)$ for the finite extension closure of the nonnegative shifts of coherent sheaves. 
Equivalently,
it consists of objects whose standard cohomology is concentrated in non-positive degrees.

For an interval $I$ of width less than one,
by \cite[Lemma~3.4 and 4.3]{Bri},
$\cP(I)$ is quasi-abelian and its strict exact sequences correspond precisely to distinguished triangles whose three terms lie in $\cP(I)$.
In any quasi-abelian category,
strict monomorphisms and epimorphisms are closed under composition,
and strict epimorphisms remain strict under pullback.

The following lemma collects the phase inequalities used throughout the paper.

\begin{lemma}\label{prelim:phase-inequalities}
For a short exact sequence
\[
0\longrightarrow F\longrightarrow E
\longrightarrow G\longrightarrow0
\]
in a heart $\HH_d$,
or a strict exact sequence in $\cP(I)$ with $|I|<1$,
the nonzero terms satisfy
\begin{equation}\label{eq:subquot-phase}
\phip(F)\le\phip(E),\qquad \phim(G)\ge\phim(E).
\end{equation}
For any triangle $F\to E\to G\stackrel{+1}{\to}$, 
one also has
$$
\phip(E)\leq\max\left\{ \phip(F),\phip(G) \right\},\qquad \phim(E)\geq\min\left\{\phim(F),\phim(G) \right\}.
$$
Moreover,
for arbitrary nonzero $E_1,E_2$,
$\Hom(E_1,E_2)=0$ if $\phim(E_1)>\phip(E_2)$.
\end{lemma}

\begin{proof}
The subobject and quotient inequalities are \cite[Lemma~3.4]{Bri},
applied to an interval of width at most one,
and the remaining parts are standard results.
\end{proof}

\paragraph{\textbf{\textit{Stability functions and hearts.}}}
A stability function on a heart $\HH$ is a homomorphism $Z:K_0(\HH)\to\C$ such that
\[
Z(E)\in
\{re^{i\pi\phi}:r>0,\ 0<\phi\leq1\}
\qquad\text{for every $0\neq E\in\HH$}.
\]
It defines semistability by $\phi_Z(F)\leq\phi_Z(E)$ for every nonzero subobject $F\subseteq E$, 
and stability by strict inequality for proper nonzero subobjects.

A bounded heart equipped with a stability function having the Harder-Narasimhan property determines a unique pre-stability condition, 
and conversely,
by \cite[Proposition~5.3]{Bri} and \cite[Lemma~12.6]{BLMNPS}. 
The support property may likewise be checked on semistable objects of the heart.

For quasi-abelian interval categories,
we use \emph{skewed stability functions}, 
whose nonzero charges lie in a rotated semi-closed half-plane
\[
\{re^{i\pi\phi}:r>0,\ a<\phi\leq a+1\}.
\]
Semistability and HN filtrations are then defined similarly by \cite[Definition~4.4]{Bri}.

A central charge on a finite-rank lattice is called \emph{rational} if its image lies in $\Q+i\Q$, 
and \emph{discrete} if its image is a discrete subgroup of $\C$.
A rational charge on a finite-rank lattice has discrete image.

\paragraph{\textbf{\textit{Mass and distances.}}}
Let $\sigma=(Z,\cP)$ be a stability condition.
For an object $E$ with HN factors $F_1,\ldots,F_n$,
define its \emph{mass} by
\[
m_\sigma(E)=\sum_{j=1}^n|Z(F_j)|\qquad \text{and }m_\sigma(0)=0.
\]
Mass is invariant under shifts,
additive under direct sums,
and subadditive along distinguished triangles,
which means that
\[
m_\sigma(E)\leq m_\sigma(F)+m_\sigma(G)
\] for every triangle $F\to E\to G\to F[1]$ by \cite[Theorem~D]{HLRMass}.
\begin{definition}
For two slicings $\cP,\cQ$,
define their generalized distance by
\[
d(\cP,\cQ)=\sup_{0\neq E}\max\left\{|\phi_{\cP}^+(E)-\phi_{\cQ}^+(E)|,\,|\phi_{\cP}^-(E)-\phi_{\cQ}^-(E)|\right\}.
\]
By \cite[Lemma~6.1]{Bri},
for $\eta\geq0$,
$d(\cP,\cQ)\leq\eta$ if and only if $\cQ(\phi)\subseteq\cP([\phi-\eta,\phi+\eta])$ for all $\phi\in\R$.
And for stability conditions $\sigma=(Z,\cP)$ and $\rho=(W,\cQ)$,
the generalized metric
\[
d(\sigma,\rho)=\max\left\{d(\cP,\cQ),\,\sup_{0\neq E}\left|\log\frac{m_\sigma(E)}{m_\rho(E)}\right|\right\}
\in[0,+\infty]
\]
induces the Bridgeland topology \cite[Proposition~8.1]{Bri}.
\end{definition}

The following standard result will be used in the base-change and moduli arguments.

\begin{lemma}
\label{prelim:topology}
If $\sigma_1=(Z_1,\cP_1)$ and $\sigma_2=(Z_2,\cP_2)$ are two stability conditions with $Z_1=Z_2$ and $d(\cP_1,\cP_2)<1$,
then $\sigma_1=\sigma_2$.
Moreover,
if $\sigma,\rho$ lie in the same connected component,
there exists $C>0$ such that
\[
C^{-1}m_\sigma(E)\leq m_\rho(E)\leq Cm_\sigma(E)
\qquad\text{for every $E\in\mathcal{T}$.}
\]
\end{lemma}
\begin{proof}
The first assertion is \cite[Lemma~6.4]{Bri}.
For the second,
since the component is path-connected,
we choose a path $\sigma_t$ from $\sigma_0=\sigma$ to $\sigma_1=\rho$.
By continuity in the generalized metric,
there is a partition $0=t_0<\cdots<t_n=1$ such that
\[
d_{\mathrm{Stab}}(\sigma_{t_{i-1}},\sigma_{t_i})<1 \qquad(1\leq i\leq n).
\]
Hence,
for every $E\neq0$,
\[
\left|\log\frac{m_\rho(E)}{m_\sigma(E)}\right|\leq\sum_{i=1}^n\left|\log\frac{m_{\sigma_{t_i}}(E)}{m_{\sigma_{t_{i-1}}}(E)}\right|<n.
\]
The assertion follows with $C=e^n$.
\end{proof}

\paragraph{\textit{\textbf{Rotation and transport.}}}
We introduce the group action on stability conditions defined in \cite[Lemma~8.2]{Bri}. 
For $a\in\R$,
rotation sends $\sigma=(Z,\cP)$ to
\[
\sigma_a=(e^{-i\pi a}Z,\cP_a),
\qquad
\cP_a(\phi)=\cP(\phi+a),
\]
whose heart is $\cP((a,a+1])$.
Rotation preserves mass,
since it multiplies the central charge by a complex number of length one.
An exact autoequivalence $\Phi$ transports the central charge and slicing by
\[
Z_{\Phi}:=Z\circ \Phi_*^{-1},\qquad(\Phi\cP)(\phi)=\Phi(\cP(\phi)).
\]
When $\Phi$ induces an automorphism of the chosen lattice $\Lambda$,
this defines an action on $\Stab_{(\Lambda,v)}(\mathcal T)$.
In particular,
tensoring by a line bundle acts in this way on the full numerical stability space.

\paragraph{\textbf{\textit{Scalar extension of stability conditions}}}
Let $X$ be a smooth projective variety over a field $k$,
and let $\cA$ be an Azumaya algebra on $X$.
Set $$
\DD=\Perf(X,\cA)=\Db(\Coh(X,\cA)).$$
For a field extension $K/k$,
let $p:X_K=X\times_k K\to X$ be the projection,
and write $$
\cA_K=p^*\cA,\qquad \DD_K=\Perf(X_K,\cA_K).
$$

Let $Y\to X$ be a Brauer--Severi scheme associated with $\cA$,
with the appropriate Brauer-class convention.
By \cite[Theorem~6.2 and Corollary~6.5]{BerghSchnuerer},
$\DD$ is equivalent to a component of a strong semiorthogonal decomposition of $D^b(Y)$.
Consequently,
\cite[Theorem~12.17]{BLMNPS} applies to numerical stability conditions on $\DD$.

In particular,
every numerical Bridgeland stability condition $\rho=(V,\mathcal{S})$ on $\DD$ admits a canonical scalar extension
\[
\rho_K=(V_K,\mathcal S_K)
\]
to $\DD_K$ for every field extension $K/k$,
satisfying the support property with respect to the base-changed numerical lattice.
These extensions are compatible with towers of fields.
Moreover,
for every $E\in\DD$ and $\phi\in\R$,
one has
\[
E\in\mathcal S(\phi)
\iff E_K\in\mathcal S_K(\phi),
\qquad
V_K(E_K)=V(E),
\]
and the Harder--Narasimhan filtration of $E_K$ is obtained by base change from that of $E$.
In particular,
$m_{\rho_K}(E_K)=m_\rho(E)$.

We will also use the description of the extended cuts
\[
\mathcal S_K(\leq c)
=
\left\{
F\in\DD_K:
p_*F\in\widehat{\mathcal S}(\leq c)
\right\}.
\]
Here $p_*$ is the pushforward on quasi-coherent derived categories,
taking values in $D_{\mathrm{qc}}(X,\cA) \simeq\Ind(\DD)$,
and $\widehat{\mathcal S}(\leq c)$ is the coaisle of the corresponding Ind-extended $t$-structure,
as recalled in Section~\ref{subsec:Ind-extensions of bounded t-structures}.

\begingroup
\let\ManuscriptSection\section
\let\ManuscriptSubsection\subsection
\let\ManuscriptSubsubsection\subsubsection
\newcount\ManuscriptProtectedSection
\ManuscriptProtectedSection=0
\renewcommand{\section}[1]{%
  \advance\ManuscriptProtectedSection by 1
  \ifcase\ManuscriptProtectedSection
  \or
    \ManuscriptSection{Phase control and generated slicings}%
    \ManuscriptSubsection{#1}%
  \or
    \ManuscriptSubsection{#1}%
  \or
    \ManuscriptSubsection{#1}%
  \or
    \ManuscriptSection{Calibration and numerical support}%
    \ManuscriptSubsection{#1}%
  \or
    \ManuscriptSubsection{#1}%
    \let\subsection\ManuscriptSubsubsection
  \else
    \PackageError{heading-layout}{Unexpected protected section}{}
  \fi
}
\subsection{Ind-extensions of bounded t-structures}
\label{subsec:Ind-extensions of bounded t-structures}
We briefly recall the $\infty$-categorical results needed below.
All derived categories are understood with their usual stable $\infty$-categorical enhancements by using $\infty$-categorical localization,
whose homotopy categories are the corresponding triangulated categories. 
Recall that a \textit{stable} $\infty$-category admits finite limits and colimits,
with pullback and pushout squares coinciding. 
And if an $\infty$-category is equivalent to a small $\infty$-category,
we call it \textit{essentially small}.

For an essentially small stable $\infty$-category $\CC$,
let $\cP(\CC)=\operatorname{Fun}(\CC^{\mathrm{op}},\mathcal{S})$ be its presheaf $\infty$-category,
where $\mathcal{S}$ denotes the $\infty$-category of spaces.
Its \textit{Ind-completion} $\Ind(\CC)$ is the smallest full sub-$\infty$-category of $\cP(\CC)$ containing the representable presheaves and closed under filtered colimits.
Equivalently,
$\Ind(\CC)$ is obtained by freely joining filtered colimits to $\CC$.

The canonical functor $\CC\hookrightarrow\Ind(\CC)$ is fully faithful.
Recall that a category is \textit{idempotent complete} if every idempotent endomorphism splits.
When $\CC$ is idempotent complete,
it identifies $\CC$ with the \textit{compact} objects of $\Ind(\CC)$,
namely those objects $P$ for which $\operatorname{Map}(P,-)$ preserves filtered colimits.
For example,
if $X$ is smooth projective,
then $$
D^b(\Coh X)\simeq\operatorname{Perf}(X),\qquad \Ind(\operatorname{Perf}(X))\simeq D_{\mathrm{qc}}(X).
$$

Suppose that $\CC$ carries a bounded $t$-structure with \textit{aisle} $\CC^{\leq0}$. 
By \cite[Proposition~2.13]{AGH},
the subcategory $\Ind(\CC^{\leq0})$ is the aisle of an induced $t$-structure on $\Ind(\CC)$.
This $t$-structure may not be bounded on $\Ind(\CC)$. 
The canonical inclusion $\CC\hookrightarrow\Ind(\CC)$ is $t$-exact,
and the truncation functors commute with filtered colimits. 
In particular,
if $M\simeq\operatorname*{colim}_i M_i$ with $M_i\in\CC$,
then $H^0(M)\simeq\operatorname*{colim}_i H^0(M_i)$ since $H^{0}=\tau^{\geq 0}\tau^{\leq 0}\simeq \tau^{\leq 0}\tau^{\geq 0}$,
where the latter colimit is taken in the heart of the induced $t$-structure.

The following lemma identifies the heart of the induced $t$-structure and shows that the noetherian property ensures the subobjects and quotients of objects in the original heart remain in that heart.

\begin{lemma}\label{lem:ind-heart}
Let $\CC$ be an essentially small stable $\infty$-category with a bounded $t$-structure,
and let $\HH=\CC^{\heartsuit}$ be its heart. 
Let $\widehat{\HH}$ be the heart of the induced $t$-structure on $\widehat{\CC}=\Ind(\CC)$.
Then $\widehat{\HH}\simeq \Ind(\HH)$,
and if $\HH$ is noetherian,
every subobject and quotient in $\Ind(\HH)$ of an object belonging to $\HH$ also belongs to $\HH$.
\end{lemma}

\begin{proof}
    Since the Ind-construction preserves fully faithful functors,
    the inclusion $\HH\hookrightarrow\CC$ induces a fully faithful functor $\Ind(\HH)\hookrightarrow\Ind(\CC)$.
    By \cite[Proposition~2.13]{AGH},
    the canonical inclusion $\CC\hookrightarrow\Ind(\CC)$ is $t$-exact,
    so $\HH\subseteq\widehat{\HH}$.
    Moreover,
    the induced heart is closed under filtered colimits. 
    Thus we obtain a fully faithful functor $\Ind(\HH)\hookrightarrow\widehat{\HH}$.

    To prove essential surjectivity,
    let $E\in\widehat{\HH}$.
    Write $E\simeq\operatorname*{colim}_i E_i$ as a filtered colimit with $E_i\in\CC$.
    Since $H^0$ commutes with filtered colimits,
    we have $$ 
    E\simeq H^0(E) \simeq H^0\bigl(\operatorname*{colim}_i E_i\bigr) \simeq\operatorname*{colim}_i H^0(E_i).
    $$
    By $t$-exactness of the canonical inclusion,
    $H^0(E_i)\in\HH$ for every $i$.
    Hence $E\in\Ind(\HH)$,
    proving $\widehat{\HH}\simeq\Ind(\HH)$.

    Now suppose that $\HH$ is noetherian.
    Let $E\in\HH$ and $S\hookrightarrow E$ be a subobject in $\widehat{\HH}$.
    Write $S\simeq\operatorname*{colim}_i S_i$,
    where $S_i\in\HH$.
    Set $I_i=\im(S_i\to E)$.
    Since $\HH\hookrightarrow\widehat{\HH}$ is exact and fully faithful,
    each $I_i$ is a subobject of $E$ computed in $\HH$.
    These images form a directed family of subobjects of $E$. 
    Since $E$ is noetherian,
    there exists an index $i_0$ such that $I_i\subseteq I_{i_0}$ for every $i$.

    Filtered colimits are exact in $\Ind(\HH)$,
    so they commute with images. 
    As $S\hookrightarrow E$ is a monomorphism,
    we obtain $$  S\simeq\im\bigl(\operatorname*{colim}_i S_i\to E\bigr)\simeq \operatorname*{colim}_i I_i\simeq I_{i_0}.
    $$
    Thus $S\in\HH$.
    Finally,
    every quotient object lies in $\HH$ since $\HH$ is an abelian category and $\HH\hookrightarrow\widehat{\HH}$ is an exact embedding.
\end{proof}

\subsection{Noetherian rational cuts}\label{subsec:Noetherian}

\begin{lemma}\label{lem:noetherian}
Suppose $(Z,\cP)$ is a stability condition and $Z(K_0(\CC))\subseteq r_1\Z+i r_2\Z$ with $r_1,r_2>0$.
Then $\cP((0,1])$ is noetherian.
If $Z(K_0(\CC))\subseteq \frac{1}{q}(\Z+i\Z)$ for some integer $q$,
the set of $d\in\R$ such that $\cP((d,d+1])$ is noetherian is dense in $\R$.
\end{lemma}

\begin{proof}
    Set $\HH=\cP((0,1])$.
    Fix a nonzero object $E\in\HH$ and consider an ascending chain $E_1\subseteq E_2\subseteq\cdots\subseteq E$ in $\HH$;
    we need to prove that it eventually stabilizes.
    Since $\Im Z(A)\geq0$ for every $A\in\HH$,
    the sequence $\Im Z(E_j)$ is nondecreasing,
    bounded above by $\Im Z(E)$,
    and takes values in $r_2\Z_{\geq0}$.
    Hence it is constant for all $j\geq j_0$.

    Set $Q_j=E_{j+1}/E_j$ and $G_j=E/E_j$.
    For $j\geq j_0$,
    every nonzero $Q_j$ has $Z(Q_j)\in r_1\Z_{<0}$,
    so $\Re Z(Q_j)\leq-r_1$.
    The exact sequence $0\to Q_j\to G_j\to G_{j+1}\to0$ gives $\Re Z(G_{j+1})-\Re Z(G_j)\geq r_1$ whenever $Q_j\neq0$.
    Thus,
    if the chain does not eventually stabilize,
    $\Re Z(G_j)\to+\infty$,
    while $y:=\Im Z(G_j)$ remains constant for $j\geq j_0$.

    Let $p=\phim(E)\in(0,1]$.
    By
    \eqref{eq:subquot-phase}, every nonzero $G_j$ belongs to $\cP([p,1])$.
    If $p<1$,
    the phase bound implies $\Im(e^{-i\pi p}Z(G_j))\geq0$,
    and hence $\Re Z(G_j)\leq y\cot(\pi p)$.
    If $p=1$,
    every nonzero $G_j$ has negative real charge.
    In either case,
    $\Re Z(G_j)$ is bounded above,
    contradicting the preceding divergence.
    Therefore $\HH$ is noetherian.

    For the second assertion,
    choose $d\in\R$ such that
    $$
    e^{i\pi d}=\frac{u+iv}{\sqrt{u^2+v^2}}
    $$
    for some $(u,v)\in\Z^2\setminus\{(0,0)\}$.
    Such choices of $d$ form a dense subset of $\R$.
    Define $Z_d=e^{-i\pi d}Z$ and $\cP_d(\phi)=\cP(\phi+d)$.
    Then $(Z_d,\cP_d)$ is a stability condition with $\cP_d((0,1])=\cP((d,d+1])$.

    If $Z(K_0(\CC))\subseteq q^{-1}(\Z+i\Z)$,
    then multiplication by $e^{-i\pi d}$ gives $$  Z_d(K_0(\CC))\subseteq\frac{1}{q\sqrt{u^2+v^2}}(\Z+i\Z).$$
    Thus the first assertion applies to $(Z_d,\cP_d)$,
    showing that $\cP((d,d+1])$ is noetherian.
    This completes the proof.  
\end{proof}

\section{Uniform phase control for vector bundles}\label{sec:uniform}
Throughout this section,
let $k$ be an arbitrary field and $X$ a smooth projective $k$-variety.
Fix a closed embedding $i:X\hookrightarrow\PP^n_k$ for $n\geq 1$,
write $H=c_1(\OO_{\PP^n}(1))$ and $\OO_X(m)=i^*\OO_{\PP^n}(m)$.
For each sufficiently large integer $M$,
set $a_M=M^2/2$ and let $\sigma_M=\tilde{\sigma}_{X}^{a_M,0}=(Z_M,\cP_M)$ be the induced geometric stability condition on $\Db(X)$ from \cite[Theorem~6.2]{LLLMPSZ},
with
\begin{equation}
\label{eq: formula of central charge}
Z_M(F)=-\int_{\PP^n_k}e^{-ia_MH}\ch(i_*F).    
\end{equation}
The following result based on \cite{LLLMPSZ} is used in the proof of Theorem~\ref{thm:main}.

\begin{lemma}
\label{prop:ambient}
For every finite set $T\subset\Z$ and every $\eta>0$,
there exists $M_0>0$ such that for all integers $M\geq M_0$,
$m\in T$,
and nonzero $F\in\Db(X)$,
one has $a_M>a_0$ and
\begin{equation}\label{eq:line-error}
\bigl|
\phi^\pm_{a_M,0}(F(m))
-\phi^\pm_{a_M,0}(F)
\bigr|<\eta.
\end{equation}
\end{lemma}

\begin{proof}
Choose an elliptic curve $E$ over $k$;
in particular, it admits a $k$-rational point $o\in E(k)$.
Let $p_j:E^n\to E$ be the projections,
and set $$
H_E:=\sum_{j=1}^n p_j^*c_1(\OO_E(o)),$$
which coincides with \cite[(4.1)]{LLLMPSZ}.
Write $\cP^{E^n}_{a,b}$ for the continuous family of slicings constructed in \cite[Proposition~4.5]{LLLMPSZ}.

For each positive integer $M$,
let $\pi_M:E^n\to E^n$ be the multiplication isogeny sending $z$ to $Mz$,
which is finite and faithfully flat.
\cite[(4.5)]{LLLMPSZ} gives $\pi_M^\sharp\cP^{E^n}_{a,b}=\cP^{E^n}_{M^2a,M^2b}$.
Applying the preceding identity to $(a,b)=(1,\frac{2m}{M^2})$ and using \cite[Lemma~3.5(3)]{LLLMPSZ},
for every fixed $m\in\Z$,
$$
d(\cP^{E^n}_{M^2,0},\cP^{E^n}_{M^2,2m}) \leq d(\cP^{E^n}_{1,0},\cP^{E^n}_{1,\frac{2m}{M^2}}) \longrightarrow 0
$$
as $M\to\infty$,
where the convergence follows from continuity of the map constructed in \cite[Proposition~4.5]{LLLMPSZ} at $(1,0)$.

Let $g:E^n\to\PP^n_k$ be the finite quotient map of \cite[Theorem~5.6]{LLLMPSZ}.
Since $g^*H=2H_E$ numerically and $\cP^{\PP^n}_{a,b}=g_\sharp\cP^{E^n}_{2a,2b}$,
tensoring with $g^*\OO_{\PP^n}(m)$ shifts the parameter $b$ by $2m$ by \cite[(4.6)]{LLLMPSZ};
hence \cite[Lemma~3.9(1)]{LLLMPSZ} gives
\[
\cP_{a_M,0}^{\PP^n}\otimes\OO_{\PP^n}(m)=g_\sharp\cP_{M^2,2m}^{E^n}.
\]

By \cite[Theorem~6.2]{LLLMPSZ},
for sufficiently large $M$,
the induced slicing on $X$ is $\cP_M=i^\sharp\cP^{\PP^n}_{a_M,0}$.
Using the distance estimate in \cite[Lemmas~3.5(3) and 3.9(3)]{LLLMPSZ},
we obtain
$$
\begin{aligned}
d(\cP_M,\cP_M\otimes\OO_X(m))
&\quad\leq
d(\cP^{\PP^n}_{a_M,0},
  \cP^{\PP^n}_{a_M,0}\otimes\OO_{\PP^n}(m))\\
&\quad\leq d(\cP^{E^n}_{M^2,0},\cP^{E^n}_{M^2,2m})\\
&\quad\leq d(\cP^{E^n}_{1,0},\cP^{E^n}_{1,2m/M^2})
\longrightarrow 0.
\end{aligned}
$$
Finally,
the tensor action satisfies $\phi^\pm_{\cP_M\otimes\OO_X(m)}(F(m))=\phi^\pm_{\cP_M}(F)$.
Hence the left-hand side of \eqref{eq:line-error} is bounded by $d(\cP_M,\cP_M\otimes\OO_X(m))$,
independent of $F$.
Since $T$ is finite,
continuity at $(1,0)$ gives a common choice of $M_0$ for all $m\in T$.
This proves the assertion.
\end{proof}

Using truncated resolutions and the cohomological bounds in the proof of \cite[Theorem~7.8]{LLLMPSZ},
we extend the preceding phase estimate to any fixed vector bundle.

\begin{proposition}\label{lem:bundle}
For every vector bundle \(B\) on \(X\) and every \(\eps>0\),
there exists $M_0=M(B,\eps)$ such that every $M>M_0$ has the following property:
every nonzero \(\sigma^X_{a_M,0}\)-semistable object \(F\in\cP_M(\phi)\) satisfies
\begin{equation}\label{eq:bundle-error}
B\otimes F\in\cP_M([\phi-\eps,\phi+\eps]).
\end{equation}
\end{proposition}

\begin{proof}
Fix a vector bundle $B$ and
pick $m_0\gg 0$ such that $B(m_0)$ is globally generated.
Then there is a surjective morphism $L_0:=\OO_X(-m_0)^{\oplus r_0}\twoheadrightarrow B$.
Repeating this procedure,
we obtain an exact sequence of length \(N\)
\begin{equation}\label{eq:bundle-resolution} 
0\longrightarrow K_N\longrightarrow L_{N-1}\longrightarrow\cdots\longrightarrow L_0\longrightarrow B\longrightarrow0,
\end{equation}
where each \(L_j\) is a finite direct sum of \(\OO_X(-m_j)\) for some $m_j\in\Z_{>0}$,
and similarly for \(B^\vee\). 
Here we choose $N=n$,
where $n$ is the dimension of the ambient projective space fixed above.
The kernels in these resolutions are locally free,
since any surjection between vector bundles splits locally.

Since all terms in the resolution are locally free,
we know $B\otimes F$ lies in the extension closure of $L_j\otimes F[j]$ for $0\leq j\leq N-1$ and $K_N\otimes F[N]$.
By Lemma~\ref{prop:ambient},
when $M$ is large enough,
every factor $F(-m)$ satisfies $\phi^{-}_{\cP_{M}}(F(-m))>\phi-\eps$,
this means $L_j\otimes F[j]\in\cP_M([\phi-\eps,\infty))$ for all  $j$.

By shifting $F$,
we may assume $0<\phi\leq 1$. 
As established in the proof of \cite[Theorem~7.8]{LLLMPSZ},
we have $$
\begin{gathered}
    \cP_M((0,1])\subset\langle \Coh(X),\ldots,\Coh(X)[N-1] \rangle,\\
    \Coh(X)\subset\cP_M((1-N,1]).
\end{gathered}
$$
Thus $K_N\otimes F[N]$ is a finite extension of objects of the form
$$
(K_N\otimes H^{-j}(F))[N+j]\in\cP_M((1+j,N+j+1]).
$$
That shows $K_N\otimes F[N]\in\cP_M(>1)$ and then $B\otimes F\in\cP_M([\phi-\eps,\infty))$. 
Increasing $M_0$ if necessary,
the same argument applied to $B^{\vee}$ shows that for every semistable object $G\in\cP_M(\psi)$,
we know $ B^\vee\otimes G\in\cP_M([\psi-\eps,\infty)).
$

If $B\otimes F=0$,
the assertion is immediate.
Otherwise,
let $G$ be its highest Harder--Narasimhan factor,
of phase $\psi$.
If $\psi>\phi+\eps$,
the lower bound for $B^\vee$ gives $B^\vee\otimes G\in\cP_M(>\phi)$.
Hence,
by tensor adjunction,
$$
\Hom(G,B\otimes F) \simeq\Hom(B^\vee\otimes G,F)=0,
$$
contradicting the nonzero morphism from the highest HN factor $G$ to $B\otimes F$.
This completes the proof.
\end{proof}

\section{Coherent generated cuts for Azumaya modules}\label{sec:cuts}

For this section \(X\) is connected and smooth projective over an arbitrary field \(k\),
and \(\cA\) is an Azumaya algebra over $X$ with degree \(r\),
hence rank \(r^2\) over \(\OO_X\).
We use left \(\cA\)-modules and write
\[
\CC=\Db(\Coh X),\quad \DD=\Db(\Coh(X,\cA)),\quad
\widehat\CC=\Dqc(X),\quad \widehat\DD=\Dqc(X,\cA).
\]
Since $X$ is smooth projective,
we know $$
\CC=\Perf(X)=\widehat{\CC}^{c},\qquad \widehat{\CC}\simeq\Ind(\CC).
$$
We define two functors $L:\CC\to\DD$ and $U:\DD\to\CC$ as 
\[
L(F)=\cA\otimes_{\OO_X}F,\qquad U(M)=M\text{ viewed over }\OO_X.
\]
We know that $L$ is left adjoint to $U$,
both are exact functors and $U$ preserves colimits.
Moreover,
$U$ conservative,
meaning that if a morphism becomes an isomorphism under $U$,
then it is an isomorphism.
An $\cA$-module is coherent as an $\OO_X$-module if and only if it is a coherent $\cA$-module.

We first record a characteristic-independent identification of the small category.
For an object $E$ in a triangulated category,
$\operatorname{thick}(E)$ denotes the smallest thick subcategory containing $E$,
i.e.,
the closure under finite direct sums,
shifts,
cones,
and direct summands.

\begin{lemma}\label{lem:azumaya-compacts}
Let $X$ be a smooth projective scheme over an arbitrary field $k$,
and let $\cA$ be an Azumaya algebra on $X$. 
Then $$
\DD=\Perf(X,\cA)=\widehat{\DD}^{\,c}.
$$
If $G$ is a compact generator of $\widehat{\CC}$,
then $LG$ is a compact generator of $\widehat{\DD}$,
and $\DD=\operatorname{thick}(LG)$. 
Consequently,
$\widehat{\DD}\simeq\Ind(\DD)$.
\end{lemma}

\begin{proof}
    Since $X$ is smooth,
    every bounded coherent $\cA$-complex is perfect,
    see \cite[Lemma~10.19]{KuznetsovHyperplane}.
    Although the cited argument is given in a characteristic-zero setting,
    it applies over any field since $\cA$ is a matrix algebra \'etale-locally,
    and Morita equivalence reduces the assertion to the regularity of $X$.
    Thus $\DD=\Perf(X,\cA)$.

    Let $P$ be a perfect $\cA$-complex. 
    The internal derived Hom functor $R\mathcal{H}om_{\cA}(P,-)$ commutes with filtered colimits. 
    Indeed,
    this can be checked \'etale-locally,
    where $P$ is represented by a bounded complex of finitely generated projective $\cA$-modules.
    Moreover,
    since $X$ is smooth,
    derived global sections commute with filtered colimits.
    Consequently,
$$
    \begin{aligned}
        R\Hom_{\cA}\bigl(P,\operatorname{colim}_{i}M_i\bigr)&\simeq R\Gamma\bigl(X,R\mathcal{H}om_{\cA}\bigl(P,\operatorname{colim}_iM_i\bigr)\bigr)\\
        &\simeq R\Gamma\bigl(X,\operatorname{colim}_iR\mathcal{H}om_{\cA}\bigl(P,M_i\bigr)\bigr)\\
        &\simeq \operatorname{colim}_i R\Hom_{\cA}\bigl(P,M_i\bigr).
    \end{aligned}
$$
    Thus $P$ is compact,
    proving $\DD\subseteq\widehat{\DD}^{\,c}$.

    By \cite[Theorem~36.15.3]{Stacks},
    $\widehat{\CC}=D_{\mathrm{qc}}(X)$ admits a perfect compact generator $G$.
    Since $L\dashv U$ and $U$ preserves filtered colimits,
    $LG$ is compact. 
    Moreover,
    if $M\in\widehat{\DD}$ satisfies $\Hom_{\widehat{\DD}}(LG[n],M)=0$ for every $n\in\Z$,
    then adjunction gives $\Hom_{\widehat{\CC}}(G[n],UM)=0$ for every $n$.
    Since $G$ generates $\widehat{\CC}$,
    we have $UM=0$;
    conservativity of $U$ then implies $M=0$.
    Thus $LG$ is a compact generator of $\widehat{\DD}$.

    By \cite[Proposition~13.37.6]{Stacks},
    we obtain $\widehat{\DD}^{\,c}=\operatorname{thick}(LG)$.
    Since $LG\in\DD$ and $\DD$ is a thick subcategory of $\widehat{\DD}$,
    this gives
    \[ \DD\subseteq\widehat{\DD}^{\,c} =\operatorname{thick}(LG)\subseteq\DD.
    \]
    Therefore   $\DD=\Perf(X,\cA)=\widehat{\DD}^{\,c}=\operatorname{thick}(LG)$.
    Finally,
    the compact-generation theorem identifies $\widehat{\DD}$ with $\Ind(\widehat{\DD}^{\,c})=\Ind(\DD)$.
\end{proof}

Fix \(0<\eps<1\).
Applying Proposition~\ref{lem:bundle} to the underlying vector bundle of $\cA$,
choose $M$ sufficiently large and set $\sigma=(Z,\cP):=\sigma_M$. 
Then
\begin{equation}\label{eq:fixed-error}
UL(\cP(\phi))\subset\cP([\phi-\eps,\phi+\eps]).
\end{equation}
for every $\phi\in\R$.

By \cite[Theorem~6.2]{LLLMPSZ},
the stability condition satisfies the support property with respect to the finite-rank numerical quotient lattice $\Lambda_i=\operatorname{Im}(K_0(X)\to\Knum(\PP^n_k))$.
Moreover,
since $a_M=\frac{M^2}{2}$,
its central charge satisfies $Z(K_0(\CC))\subseteq \frac{1}{2^n n!}(\Z+i\Z)$.
In particular,
$Z$ has discrete image.
We keep this stability condition fixed for all cutoffs $c\in\R$.

To construct a stability condition on $\DD$, 
a natural candidate for the central charge is $W=Z\circ U$, 
where $U$ is the forgetful functor.
The main difficulty is to construct a suitable slicing compatible with $W$.
Indeed, the naive choice $U^{-1}\cP(\phi)$ fails since it need not satisfy the Harder--Narasimhan property.

An easy counter example is $X=\PP_k^{1}$, 
$V=\OO_X\oplus\OO_X(1)$ and $\cA=\mathcal{E}nd_{\OO_X}(V)$.
Let $\sigma_t=(Z_t,\cP_t)$ be the standard slope stability condition on $X$, with $Z_t(F)=-\deg(F)+it\rk(F)$ for all $t>0$.
Then $UL(\OO)=\OO(-1)\oplus\OO^{\oplus 2}\oplus\OO(1)\notin \cP_t(\frac{1}{2})$.

Moreover, the subcategories $U^{-1}\cP_t(\phi)$ do not form a slicing. 
Indeed, every object $E\in U^{-1}\cP_t(\phi)$ is of the form $V\otimes F$ and then $UE=F\oplus F(1)\in \cP_t(\phi)$ for some $F\in\CC$ by Morita equivalence.
Since $Z_t(F(1))=Z_t(F)-\rk(F)$ and $F,F(1)$ has the same phase since they both lie in $\cP_t(\phi)$. 
So we know $\rk(F)=0$.
Consider $E_0=V$, if $\{U^{-1}\cP_t(\phi)\}_{\phi\in\R}$ form a slicing, the same argument implies all HN-factors of $E_0$ are of rank zero, gives a contradiction.

This counter example suggests a loose version of phase control:
instead of forcing $UE$ to be semistable of phase $\phi$,
we seek a slicing $\cQ$ such that all HN phases of $UE$ lies in a small interval containing $\phi$. And then modify the obtained slicing $\cQ$ appropriately to get the stability condition $(W,\cR)$.
Moreover, the counter example above told us that we cannot lift the HN-filtration of $UE\in\CC$ directly back to $\DD$, instead we want to lift the $t$-structures determined by $\cP$.
Here the technique is to consider the Ind-extension of $\cP$ to $\widehat{\CC}$, pullback by $U$ to $\widehat{\DD}$ and then restrict back to $\DD$.
We then use the controlled phase errors to construct the locally finite slicing $\cQ$.

For any $c\in\R$,
the slicing $\cP$ defines a $t$-structure on $\CC$ with aisle $\cP(>c)$ and coaisle $\cP(\leq c)$;
we denote its Ind-extension to $\widehat{\CC}$ by $\widehat{\cP}(>c)=\Ind(\cP(>c))$ and $\widehat{\cP}(\leq c)=\widehat{\cP}(>c)^{\perp}$.
Let \(\widehat\cA_c\) be the aisle in \(\widehat\DD\) generated under extensions and colimits by
\[
\mathcal{S}_c=\{L(F):F\in\cP(>c)\}\subseteq\widehat{\DD},
\]
and let \(\widehat\BB_c=\widehat\cA_c^\perp\).
By \cite[Proposition~1.4.4.11]{HA} applied to $\mathcal{T}=\widehat{\DD}$ and $\mathcal{S}_c$,
one obtains a $t$-structure $(\widehat{\cA}_c,\widehat{\BB}_c)$ satisfying $$
\widehat{\cA}_c[1]\subset\widehat{\cA}_c,\qquad\Hom(\widehat{\cA}_c,\widehat{\BB}_c)=0,
$$ and for any $M\in\widehat{\DD}$,
there exists a truncation triangle $M_+\to M\to M_{-}$,
where $M_+\in\widehat{\cA}_c,M_{-}\in\widehat{\BB}_c$.

The following lemma shows the relation between Ind-extension of $t$-structure induced by slicing $\cP$ and $(\widehat{\cA}_c,\widehat{\BB}_c)$.
\begin{lemma}\label{lem:large-bands}
For every $c\in\R$,
the generated $t$-structure satisfies
\begin{equation}\label{eq:large-bands}
\widehat\BB_c=U^{-1}\widehat\cP(\le c),\qquad U(\widehat\cA_c)\subseteq\widehat\cP(>c-\eps).
\end{equation}
Moreover,
these $t$-structures are nested: $\widehat{\cA}_b\subseteq\widehat{\cA}_a$ and $\widehat{\BB}_a\subseteq\widehat{\BB}_b$ whenever $a<b$, 
and satisfy $\widehat{\cA}_{c+1}=\widehat{\cA}_c[1]$ and $\widehat{\BB}_{c+1}=\widehat{\BB}_c[1]$.
\end{lemma}

\begin{proof}
    By the universal property of the $t$-structure generated by $L(\cP(>c))$,
    an object $M\in\widehat{\DD}$ belongs to $\widehat{\BB}_c$ if and only if
    \[ \Hom_{\widehat{\DD}}(LF,M)=0\qquad \text{for every }F\in\cP(>c).
    \]
    By adjunction,
    this is equivalent to $\Hom_{\widehat{\CC}}(F,UM)=0$ for every $F\in\cP(>c)$.
    Since $\widehat{\cP}(>c)$ is the Ind-extension of $\cP(>c)$,
    the latter condition is equivalent to $UM\in\widehat{\cP}(\leq c)$.
    Consequently,
    $\widehat{\BB}_c  =U^{-1}\widehat{\cP}(\leq c)$.

    Every generator of $\widehat{\cA}_c$
    can be written as $L(F)$ for some $F\in\cP(>c)$.
    Let $F_i$ be its HN-factors with phase $\phi_i$.
    By \eqref{eq:bundle-error},
    we know $UL(F_i)\in\cP([\phi_i-\eps,\phi_i+\eps])\subseteq\cP(>c-\eps)$.
    Since $UL$ is exact and $\cP(>c-\eps)$ is closed under extension,
    we know $UL(F)\in\cP(>c-\eps)$ and denote $$\mathcal{T}_c=\left\{M\in\widehat{D}:UM\in\widehat{\cP}(>c-\eps)  \right\}.$$
    It is easy to check $\mathcal{T}_c$ is closed under extensions and colimits,
    so $\widehat{\cA}_c\subseteq \mathcal{T}_c$ by definition of $\widehat{\cA}_c$.

    Finally,
    $\cP(>b)\subseteq\cP(>a)$ for $a<b$ gives $\widehat{\cA}_b\subseteq\widehat{\cA}_a$,
    and taking right orthogonals yields $\widehat{\BB}_a\subseteq\widehat{\BB}_b$.
    Since $\cP(>c+1)=\cP(>c)[1]$ and $L$ is exact,
    we have $\widehat{\cA}_{c+1}=\widehat{\cA}_c[1]$.
    Taking right orthogonal complements gives $\widehat{\BB}_{c+1}=\widehat{\BB}_c[1]$.
\end{proof}

The preceding lemma controls the generated $t$-structures through the forgetful functor. 
We now show that their truncations preserve bounded coherent objects,
using the noetherian ambient hearts from Section~\ref{subsec:Ind-extensions of bounded t-structures}--\ref{subsec:Noetherian}.

\begin{proposition}\label{prop:coherence}
For every $c\in\R$,
the generated $t$-structure $(\widehat{\cA}_c,\widehat{\BB}_c)$ restricts to a bounded $t$-structure $(\cA_c,\BB_c)$ on $\DD$,
where $\cA_c=\widehat{\cA}_c\cap\DD$ and $\BB_c=\widehat{\BB}_c\cap\DD$.
\end{proposition}

\begin{proof}
Let \(M\in\DD\),
and take its truncation triangle
\[
M_+\longrightarrow M\longrightarrow M_-\longrightarrow M_+[1]
\]
in $\widehat{\DD}$,
where $M_+\in\widehat{\cA}_c$ and $M_-\in\widehat{\BB}_c$.
Applying forgetful functor $U$,
we obtain a triangle
$$
UM_+\longrightarrow UM\longrightarrow UM_-\longrightarrow UM_+[1]
$$
and by Lemma~\ref{lem:large-bands},
we know $UM_+\in\widehat{\cP}(>c-\eps)$ and $UM_-\in\widehat{\cP}(\leq c)$.
Choose \(d\in\R\) from the dense subset of $\R$ in Lemma~\ref{lem:noetherian} such that $c-1<d<c-\eps$.
The heart \(\HH_d=\cP((d,d+1])\) is noetherian,
and the induced $t$-structure on $\widehat{\CC}$ has heart $\widehat{\HH}_d\simeq\Ind(\HH_d)$.
Consider the bounded $t$-structure $(\CC_d^{\leq 0},\CC_d^{\geq 0})$ induced by $\HH_d$,
more precisely,
$\CC_d^{\leq 0}=\cP(>d)$,
$\CC_d^{\geq 0}=\cP(\leq d+1)$.
And denote $\widehat{\CC}_d^{\leq 0}=\widehat{\cP}(>d)$ and $\widehat{\CC}_d^{\geq 0}=\widehat{\cP}(\leq d+1)$.

By definition above, we know \[
UM_+\in\widehat\CC_d^{\le0},\qquad
UM_-\in\widehat\CC_d^{\ge0}.
\]
Let $H_d^j:\widehat{\CC}\to\widehat{\HH}_d$ be the cohomology functor with respect to the heart $\HH_d$. 
Since $H_d^{-1}(UM_-)=H_d^1(UM_+)=0$,
it induces a short exact sequence 
\begin{equation}\label{eq:degree-zero}
0\longrightarrow H_d^0(UM_+)\longrightarrow H_d^0(UM)\longrightarrow H_d^0(UM_-)\longrightarrow0.
\end{equation}
The middle term lies in \(\HH_d\) since $M\in \DD$ implies $UM\in\CC$,
so Lemma~\ref{lem:ind-heart} implies that $H_d^0(UM_+)$ and $H_d^0(UM_-)$ lies in $\HH_d$ as well,
in particular,
they are compact.

Since $UM_+\in\widehat{\CC}_d^{\leq 0}$ and $UM_-\in\widehat{\CC}_d^{\geq 0}$,
by the universal properties of truncation,
we know
\[
\tau_d^{\le-1}(UM_+)\simeq\tau_d^{\le-1}(UM),\qquad\tau_d^{\ge1}(UM_-)\simeq\tau_d^{\ge1}(UM).
\]
Consequently there are two canonical triangles
\begin{align*}
\tau_d^{\le-1}(UM)&\longrightarrow UM_+
\longrightarrow H_d^0(UM_+),\\
H_d^0(UM_-)&\longrightarrow UM_-\longrightarrow\tau_d^{\ge1}(UM).
\end{align*}

Since the inclusion $j:\CC\hookrightarrow\widehat{\CC}$
is $t$-exact for the induced $t$-structure,
we know $$
j(\tau_{d,\CC}^{\leq -1}(UM))\simeq \tau_{d,\widehat{\CC}}^{\leq -1}(j(UM)),\qquad j(\tau_{d,\CC}^{\geq 1}(UM))\simeq \tau_{d,\widehat{\CC}}^{\geq 1}(j(UM)).
$$

We know $H^0_d(UM_+),H^0_d(UM_-)\in\HH_d\subseteq\CC$,
and $\tau_d^{\leq-1}(UM),\tau_d^{\geq1}(UM)\in\CC$,
so $UM_+,UM_-\in\CC$,
which means that $UM_+$ and $UM_-$ are bounded coherent,
then $M_+,M_-\in\DD$ are also bounded coherent.
Here smoothness of $X$ ensures that these bounded coherent $\mathcal A$-complexes are perfect, hence compact, as in Lemma~\ref{lem:azumaya-compacts}.

Finally,
we need to prove for any $M\in\DD$,
$M\in\BB_c$ for $c$ large enough and $M\in\cA_c$ for $c$ small enough.
Since $UM\in\CC$,
it has a finite HN-filtration,
so there exists $c_+\in\R$ such that $UM\in\cP(\leq c)$ for $c\geq c_+$,
then $M\in\widehat{\BB}_c\cap\DD=\BB_c$ for all sufficiently large $c$.
To prove the other bound,
set
\[
\mathcal T=\bigcup_{c\in\R}\cA_c\subset\DD.
\]
This is a thick triangulated subcategory of $\DD$,
which follows from the nesting and shift compatibility of the aisles $\cA_c$ for all $c$,
together with their closure under extensions and retracts. 
Picking a compact generator $G$ of $\widehat{\CC}$,
we know $\DD=\operatorname{thick}(LG)$ by Lemma~\ref{lem:azumaya-compacts}. 
Since $G\in\CC$,
it has finite HN-filtration,
so pick $c<\phi_{\cP}^{-}(G)$,
then $G\in\cP(>c)$,
so $LG\in\widehat{\cA}_c\cap\DD=\cA_c\subset\mathcal{T}$. 
Because $\mathcal{T}$ is thick,
$\DD=\operatorname{thick}(LG)\subseteq\mathcal{T}$,
so $\mathcal{T}=\DD$.
As a result,
$M\in\cA_c$ for all sufficiently small $c$.

For any $c\in\R$ and $M\in\cA_{c_-}$,
$M\in\BB_{c_+}$ for some $c_-,c_+$.
Pick $n$ large enough such that $c-n<c_-$ and $c+n>c_+$,
so $M\in\cA_{c}[-n]\cap\BB_{c}[n]$.
So $(\cA_c,\BB_c)$ is bounded,
that completes the proof.
\end{proof}

\section{Finite transitions and a slicing}\label{sec:transitions}
We now show that the bounded $t$-structures constructed in Section~\ref{sec:cuts} arise from a locally finite slicing on $\DD$.
Recall that the fixed stability condition $\sigma=(Z,\cP)$ satisfies the support property with respect to a finite-rank lattice by \cite[Theorem~6.2]{LLLMPSZ}.
By Lemma~\ref{lem:large-bands} and Proposition~\ref{prop:coherence},
the bounded $t$-structures $(\cA_c,\BB_c)$ satisfy
\begin{equation}\label{eq:cut-phase-range}
U(\cA_a\cap\BB_b)
\subseteq\cP((a-\eps,b])
\qquad(a<b).
\end{equation}

\begin{proposition}
\label{prop:finite-transitions}
For each $\phi\in\R$,
define
\begin{equation}\label{eq:twisted-slicing}
\cQ(\phi)
=
\left(\bigcap_{a<\phi}\cA_a\right)
\cap
\left(\bigcap_{b>\phi}\BB_b\right).
\end{equation}
Then $\cQ$ is a locally finite slicing on $\DD$ satisfying
\begin{equation}\label{eq:slicing-phase-range}
U(\cQ(\phi))
\subseteq\cP([\phi-\eps,\phi])
\qquad(\phi\in\R).
\end{equation}
Moreover,
$\cQ(I)$ has finite length for every interval $I$ of width less than $1-\eps$.
The $t$-structures induced by $\cQ$ coincide with the $t$-structures constructed in Section~\ref{sec:cuts};
namely,
$\cQ(>c)=\cA_c$ and $\cQ(\leq c)=\BB_c$ for every $c\in\R$.
\end{proposition}

\begin{proof}
\textit{Step 1:} Choose a lattice norm and a constant \(C_\sigma\) as in \eqref{eq:support},
and set \(\lambda:=\inf_{\gamma\in\Lambda\setminus\{0\}}\|\gamma\|>0\).
Since the central charge factors through $\Lambda$,
every nonzero $\sigma$-semistable object $G$ satisfies $|Z(G)|\geq \frac{\lambda}{C_{\sigma}}$.

For a closed interval \(J=[u,v]\) of width \(w<1\),
put \(\theta=(u+v)/2\) and define the additive function $\rho_J(F)=\Re(e^{-i\pi\theta}Z(F))$ on $K_0(\CC)$.
If \(0\ne F\in\cP(J)\),
its HN factors have phases in $J$,
and hence
\begin{equation}\label{eq:positive-weight}
\begin{aligned}
\rho_J(F) &\geq \cos\left(\frac{\pi w}{2}\right)m_\sigma(F)
\geq \frac{\lambda}{C_\sigma} \cos\left(\frac{\pi w}{2}\right)=:\mu_J>0.
\end{aligned}
\end{equation}
In particular,
$\rho_J$ is additive on distinguished triangles and has a uniform positive lower bound on the nonzero objects of $\cP(J)$.

\smallskip
\noindent\emph{Step 2:}
Fix $E\in\DD$. 
For each $c\in\R$,
let $H_c=H_c(E)$ and $K_c=K_c(E)$ be the terms of its canonical truncation triangle $H_c\to E\to K_c$,
with $H_c\in \cA_c$ and $K_c\in\BB_c$. 
For $a<b$,
by the universal property of the truncation functor,
there is a canonical map \(H_b\to H_a\);
denote its cone by $E_{a,b}$.

We claim that $E_{a,b}\in\cA_a\cap \BB_b$.
Indeed,
$H_a,H_b[1]\in\cA_a$,
so extension closure gives $E_{a,b}\in \cA_a$. 
On the other hand,
the octahedral axiom gives a triangle $K_a[-1]\to E_{a,b}\to K_b$.
Since $K_a[-1],K_b\in\BB_b$,
it follows that $E_{a,b}\in\BB_b$.

Now fix $a<b$ with $b-a+\eps<1$,
and put $J=[a-\eps,b]$.
For any subdivision
\(a=c_0<\cdots<c_n=b\), the octahedral axiom gives $$
[E_{a,b}]=\sum_{j=1}^n[E_{c_{j-1},c_j}]\quad\text{in } K_0(\DD).
$$
By \eqref{eq:cut-phase-range},
every $U(E_{c_{j-1},c_j})$ belongs to $\cP(J)$.
Moreover,
conservativity of $U$ and \eqref{eq:positive-weight} give $\rho_J(U(E_{c_{j-1},c_j}))\geq\mu_J$ whenever $E_{c_{j-1},c_j}\ne0$. 
Thus \begin{equation}\label{eq:transition-bound}
\#\{j:E_{c_{j-1},c_j}\ne0\}\leq\frac{\rho_J(U(E_{a,b}))}{\mu_J}<\infty.
\end{equation}
In particular,
the number of nontrivial cutoff changes in any subdivision of $[a,b]$ is uniformly bounded.

Call two cutoffs $s<t$ \textit{equivalent} if the canonical map $f_{t,s}:H_t\to H_s$ is an isomorphism.
Their equivalence classes are convex intervals.
Indeed,
if $s<u<t$ and $f_{t,s}:H_t\to H_s$ is an isomorphism,
then the truncation triangle $H_s\to E\to K_s$ at $s$ is also a truncation triangle for $(\cA_u,\BB_u)$. 
Uniqueness of truncations therefore identifies all three cutoffs.

If $\left\{H_c\right\}_{c\in [a,b]}$ has infinitely many classes, 
we can put arbitrarily many points $a\leq t_0<\cdots< t_m\leq b$ such that $E_{t_{j-1},t_j}\neq 0$ for each $j$,
contradicting with \eqref{eq:transition-bound}, 
hence there are only finitely many classes of $H_c$ whenever $b-a+\eps<1$. 

Since $(\cA_0,\BB_0)$ is bounded by Section~\ref{sec:cuts}, 
there exists an integer $n$ large enough such that $E\in \cA_{-n}$,
so $H_c\simeq E$ for $c\ll 0$,
and similarly $E\in \BB_m$ for $m$ large enough,
so $H_c\simeq 0$ when $c\gg 0$.

\smallskip
\noindent\emph{Step 3:}
By the above arguments,
we know for any fixed object $E$,
$H_c(E)$ varies finitely many times,
and we denote by $\phi_1<\phi_2<\ldots<\phi_r$ to be the points where $H_c(E)$ changes.
Fix a point $\phi=\phi_i$,
choose $\delta>0$ such that $\phi_i\in(\phi_{i-1}+\delta,\phi_{i+1}-\delta)$.
When $s\in\left(\phi-\delta,\phi\right)$,
denote $H_s(E)=:H_{\phi^{-}}$ and for $s\in\left(\phi,\phi+\delta\right)$,
denote $H_s(E)=:H_{\phi^{+}}$ and $C_{\phi}:=\Cone(H_{\phi^+}\to H_{\phi^-})$. 
For any $\phi-\delta<s<\phi<t<\phi+\delta$,
$C_{\phi}\simeq E_{s,t}$.
Consequently,
$C_{\phi}\in \cA_s\cap \BB_t$ for such $s,t$,
means $C_{\phi}\in\cQ(\phi)$ when $s\uparrow\phi$ and $t\downarrow\phi$.

Choosing cutoffs $$
s_0<\phi_1<s_1<\phi_2<\cdots<\phi_r<s_r
$$
in the intervening constant regions,
we know there exist canonical maps $0=H_{s_r}\to H_{s_{r-1}}\to\cdots\to H_{s_0}=E$, 
the successive cones $C_j=\Cone(H_{s_j}\to H_{s_{j-1}})$ lies in $\cQ(\phi_j)$,
it forms a finite $\cQ$-HN filtration of $E\in \DD$.
Thus every nonzero object of $\DD$ admits a finite $\cQ$-HN filtration.

Moreover,
since $\cA_a,\BB_b$ is additive,
extension-closed subcategories for any $a,b$,
so $\cQ(\phi)$ are additive and extension-closed subcategories.
What's more,
by Lemma~\ref{lem:large-bands},
$\cA_c[1]=\cA_{c+1}$ and $\BB_{c}[1]=\BB_{c+1}$,
so $\cQ(\phi+1)=\cQ(\phi)[1]$.
For $\phi>\psi$,
choosing $\psi<c<\phi$ gives $\cQ(\phi)\subseteq\cA_c$ and $\cQ(\psi)\subseteq\BB_c$, 
hence $\Hom(\cQ(\phi),\cQ(\psi))=0$.
These properties prove that $\cQ$ is a slicing.

For $0\ne F\in\cQ(\phi)$ and arbitrary $a<\phi<b$,
we have $UF\in\cP((a-\eps,b])$ by \eqref{eq:cut-phase-range}.
Taking $a\uparrow\phi$ and $b\downarrow\phi$ in the inequalities gives $UF\in\cP([\phi-\eps,\phi])$.
This proves \eqref{eq:slicing-phase-range}.

\smallskip
\noindent\emph{Step 4:}
Let $I$ be an interval of width less than $1-\eps$,
and put $J=[\inf I-\eps,\sup I]$. 
Then $J$ has width less than one,
and \eqref{eq:slicing-phase-range} together with the exactness of $U$ gives $U(\cQ(I))\subseteq\cP(J)$.

By \cite[Lemma~4.3]{Bri},
the category $\cQ(I)$ is quasi-abelian.
Every nonzero object $F\in\cQ(I)$ satisfies $\rho_J(UF)\geq\mu_J$,
since $U$ is conservative and \eqref{eq:positive-weight}.
Thus any finite strict filtration of $E\in\cQ(I)$ has at most $\frac{\rho_J(UE)}{\mu_J}$ nonzero factors.
This bounds all strict subobject and quotient chains,
proving that $\cQ(I)$ has finite length.
In particular,
$\cQ$ is locally finite.

Finally,
Lemma~\ref{lem:large-bands} gives $\BB_c=U^{-1}\cP(\leq c)$.
We first observe that $ \widehat{\cP}(\leq c)=\bigcap_{b>c}\widehat{\cP}(\leq b).
$ Using $\widehat{\BB}_c =U^{-1}\widehat{\cP}(\leq c)$ and restricting to $\DD$, we obtain
\[
\BB_c=\bigcap_{b>c}\BB_b.
\]

By definition of $\cQ$ and this right-continuity,
$\cQ(>c)\subseteq\cA_c$ and $\cQ(\leq c)\subseteq\BB_c$.
Hence the $\cQ$-HN cutoff triangle at $c$ is also a truncation triangle for $(\cA_c,\BB_c)$.
Uniqueness of truncations gives $\cQ(>c)=\cA_c$ and $\cQ(\leq c)=\BB_c$,
as required.
\end{proof}

\begin{remark}\label{rem:abstract-slicing}
The proof of Proposition~\ref{prop:finite-transitions} uses only the following hypotheses: 
\begin{itemize}
    \item An exact conservative functor $U:\DD\to\CC$.
    \item A stability condition $(Z,\cP)$ on $\CC$ satisfying the support property on a finite-rank lattice.
    \item A family of bounded $t$-structures $(\cA_c,\BB_c)_{c\in\R}$ such that for $a<b$, $\cA_b\subseteq\cA_a$,
    $\BB_a\subseteq\BB_b$,
    $\cA_{c+1}=\cA_c[1]$,
    and $U(\cA_a\cap\BB_b)\subseteq\cP((a-\eps,b])$ for some $0\leq\eps<1$.
\end{itemize}

Under these hypotheses,
Proposition~\ref{prop:finite-transitions} defines a locally finite slicing $\cQ$ satisfying \eqref{eq:slicing-phase-range}, 
and $\cQ(I)$ has finite length whenever $I$ has width less than $1-\eps$.
If,
in addition,
$\BB_c=\bigcap_{b>c}\BB_b$,
its induced cutoffs coincide with the original $t$-structures.

In particular, the construction of $\cQ$ does not intrinsically depend on the smoothness of $X$. 
This raises the question of whether the same method applies to singular projective schemes.
The main difficulty is to construct the required bounded $t$-structures on the appropriate category, 
rather than recover a slicing from such $t$-structure.

Indeed, the proof of Proposition~\ref{prop:coherence} uses smoothness to identify $$
\Db(\Coh(X,\cA))=\Perf(X,\cA)=\mathrm{D}_{\mathrm{qc}}(X,\cA)^c.
$$
For singular $X$, bounded coherent complexes need not be perfect, and the generated $t$-structures on $\mathrm{D}_{\mathrm{qc}}(X,\mathcal A)$ need not restrict to $\operatorname{Perf}(X,\mathcal A)$. 
Nevertheless, once the suitable bounded cuts is formulated, the argument of Proposition~\ref{prop:finite-transitions} still applies.
\end{remark}

\begin{corollary}
\label{cor:induced-charge}
Let $W=Z\circ U:K_0(\DD)\to\C$.
Then $W(K_0(\DD))$ is discrete.
Moreover,
for every nonzero $E\in\cQ(\phi)$,
one has $W(E)\ne0$ and $\phi_W(E)\in[\phi-\eps,\phi]$.
\end{corollary}

\begin{proof}
By Proposition~\ref{prop:finite-transitions},
$UE$ is nonzero and belongs to $\cP([\phi-\eps,\phi])$.
Since $\eps<1$,
the central charges of nonzero HN factors of $UE$ lie in a sector of angle less than $\pi$.
Their sum $W(E)=Z(UE)$ is therefore nonzero and has phase $\phi_W(E)\in[\phi-\eps,\phi]$.

Finally,
the chosen ambient stability condition $\sigma=(Z,\cP)$ has discrete central charge image,
since its parameter is rational.
Thus $W(K_0(\DD))\subseteq Z(K_0(\CC))$ is discrete.
\end{proof}

\section{Phase calibration for a discrete charge}\label{sec:calibration}
In Section~\ref{sec:transitions},
we obtain a slicing $\cQ$ on $\DD=\Db(\Coh(X,\cA))$,
a group homomorphism $W:K_0(\DD)\to\C$ with discrete image and a real number $\delta>0$ small enough such that every $0\neq E\in \cQ(\phi)$ satisfies $W(E)\neq 0$ and 
\begin{equation}\label{eq: argument of W}
\phi-\delta<\frac{1}{\pi}\arg W(E)<\phi+\delta.    
\end{equation}
Here the argument is understood as a real lift chosen in the appropriate interval.
We now adjust $\cQ$ to make it compatible with the induced central charge $W=Z\circ U$.
The construction is an adaptation of \cite[Section~7]{Bri}.
Although \cite[Theorem~7.1]{Bri} cannot apply directly,
\eqref{eq: argument of W} and Proposition~\ref{prop:finite-transitions} give the necessary inputs for the argument in this section.

\begin{theorem}\label{thm:calibration}
Let \(0<\delta<\frac{1}{16}\),
and $\cQ,W$ be as above. 
Then there exists a locally finite pre-stability condition \((W,\cR)\) on $\DD $ satisfying
\[
\cR(\psi)\subset\cQ((\psi-\delta,\psi+\delta))\, \text{ for all } \psi\in\R, \qquad d(\cQ,\cR)\le\delta.
\]
\end{theorem}

\begin{proof}
We adapt the proof technique of \cite[Section~7]{Bri}.

\smallskip
\noindent\emph{Step 1:}
Call $\mathcal E=\cQ((a,b))$ \emph{thin} if $0<b-a<1-2\delta$.
By \cite[Lemma~4.3]{Bri},
$\mathcal E$ is a quasi-abelian category whose exact sequences are given by distinguished triangles.

Let $0\neq F\in\mathcal E$,
let $F_j\in\cQ(\phi_j)$ be its $\cQ$-HN factors,
where $a<\phi_j<b$.
By \eqref{eq: argument of W},
we know $\frac{1}{\pi}\arg W(F_j)\in(a-\delta,b+\delta)$ with width $w=b-a+2\delta<1$.
Hence $W(F)\neq 0$ and there exists a unique $\psi(F)\in (a-\delta,b+\delta)$ such that $$W(F)=|W(F)|e^{i\pi \psi(F)}.$$
This defines a skewed stability function on $\mathcal E$,
as in \cite[Definition~7.2]{Bri}.

We claim that every thin interval has finite length.
Set \[
\theta=\frac{a+b}{2},\qquad
\kappa=\cos\left(\frac{\pi(b-a+2\delta)}{2}\right)>0,
\] and define $\rho(F)=\Re(e^{-i\pi\theta}W(F))$.
Since $W(K_0(\DD))$ is discrete by Corollary~\ref{cor:induced-charge}, 
there exists $m_W>0$ such that $|z|\geq m_W$ for every nonzero $z\in W(K_0(\DD))$. 
Hence \[
\rho(F)\geq\kappa|W(F)|\geq\kappa m_W>0.
\]
Since $\rho$ is additive in strict exact sequences,
every strict filtration of a fixed object $F$ has at most $\frac{\rho(F)}{\kappa m_W}$ nonzero factors.
Thus $\mathcal E$ has finite length.

\smallskip
\noindent\emph{Step 2:}
A nonzero object $F\in\mathcal E$ is called \textit{$W$-semistable} if $\psi(A)\leq\psi(F)$ for every nonzero strict subobject $A\hookrightarrow F$ in $\mathcal E$.
Following \cite[Definition~7.4]{Bri},
we say $F$ is \emph{enveloped} by $\mathcal E$ if $$
a+\delta\leq\psi(F)\leq b-\delta.
$$

The proof of \cite[Lemma~7.3]{Bri} applies in our setting,
if $F$ is $W$-semistable in a thin interval,
then 
\begin{equation}\label{eq: phase of F}
F\in\cQ((\psi(F)-\delta,\psi(F)+\delta)).  \end{equation}
Likewise,
the argument of \cite[Lemma~7.5]{Bri} shows that $W$-semistability of an enveloped object is independent of the chosen enveloping interval.
In fact,
the proofs of these two lemmas use only the slicing axioms and \eqref{eq: argument of W}.

For each $\psi\in\R$,
let $\cR(\psi)$ consist of zero and the objects that are $W$-semistable of phase $\psi$ in some thin enveloping interval .
By \eqref{eq: phase of F},
\begin{equation}\label{eq:calibration-confinement}
\cR(\psi)\subseteq
\cQ((\psi-\delta,\psi+\delta)).
\end{equation}
The proof of \cite[Lemma~7.6]{Bri},
using kernels and images in a common $\cQ$-heart,
then gives
\begin{equation}\label{eq:calibration-hom}
\Hom_{\DD}(\cR(\psi_1),\cR(\psi_2))=0
\end{equation}
for any $\psi_1>\psi_2$.
The definition of $\cR(\psi)$ also implies $\cR(\psi+1)=\cR(\psi)[1]$.

Finally,
the proof of \cite[Lemma~7.7]{Bri} applies, 
since its argument uses the phase estimate \eqref{eq: argument of W},
the Hom-vanishing \eqref{eq:calibration-hom}, 
and the finite-length property of thin intervals established in Step 1.
So every object of $\cQ((a+2\delta,b-4\delta))$ has a finite HN filtration in $\mathcal{E}$,
whose factors are $W$-semistable and enveloped by $\mathcal{E}$,
hence belong to the subcategories $\cR(\psi)$ for some $\psi$.
We refer to this as the \textit{buffered HN construction}.

\smallskip
\noindent\emph{Step 3:}
For $t\in\R$,
define
\[
\mathcal U_t=
\langle\cR(\psi):\psi>t\rangle_{\mathrm{ext}},
\qquad
\mathcal V_t=
\langle\cR(\psi):\psi\leq t\rangle_{\mathrm{ext}},
\]
By \eqref{eq:calibration-hom},
$\Hom(\mathcal U_t,\mathcal V_t)=0$.
Shift compatibility of $\cR$ gives $\mathcal{U}_t[1]\subseteq \mathcal{U}_t$.

For $0\ne F\in \cQ(s)$, 
apply the buffered HN construction in the thin interval $\cQ((s-3\delta,s+5\delta))$.
Let $G_1,\ldots, G_{\ell}$ be the resulting HN factors,
with $G_j\in\cR(\psi_j)$ and $\psi_1>\cdots>\psi_{\ell}$.
The first and last factors give nonzero morphisms $$
\Hom(G_1,F)\neq 0,\qquad \Hom(F,G_{\ell})\ne 0.
$$
By \eqref{eq:calibration-confinement},
$G_j\in\cQ((\psi_j-\delta, \psi_j+\delta))$.
Since $F\in\cQ(s)$,
we know $\psi_1-\delta<s$ and $\psi_{\ell}+\delta>s$,
so we obtain $s-\delta<\psi_j<s+\delta$ for all $1\leq j\leq \ell$.

Consequently,
taking finite extension yields
\begin{equation}\label{eq:calibration-cut-bounds}
\cQ([t+\delta,\infty))\subseteq\mathcal U_t,
\qquad
\cQ((-\infty,t-\delta])\subseteq\mathcal V_t.
\end{equation}

We now construct cut triangles,
following the proof of \cite[Theorem~7.1]{Bri}.
Consider the $\cQ$-HN filtration,
we know any $M\in\DD$ belongs to some $H*C*L$,
where
\[
\begin{aligned}
H&\in\cQ([t+\delta,\infty)),\\
C&\in\cQ((t-\delta,t+\delta)),\\
L&\in\cQ((-\infty,t-\delta]).
\end{aligned}
\]
By \eqref{eq:calibration-cut-bounds},
we know $H\in \mathcal{U}_t$ and $L\in\mathcal{V}_t$.
For $C\in\cQ((t-\delta,t+\delta))$,
apply the buffered HN construction in $\cQ((t-3\delta,t+5\delta))$,
we obtain a finite HN-filtration of $C$,
and partition these factors by phases $>t$ or $\leq t$,
we know $C\in\mathcal{U}_t*\mathcal{V}_t$,
so $M\in \mathcal{U}_t *\mathcal{V}_t$,
means $M$ admits a triangle $M_t^+\to M\to M_t^-$ with $M_t^+\in \mathcal{U}_t$ and $M_t^-\in\mathcal{V}_t$.

To prove $(\mathcal{U}_t,\mathcal{V}_t)$ is a $t$-structure, 
it remains to show $\mathcal{U}_t^{\perp}\subseteq \mathcal{V}_t$,
for $M\in \mathcal{U}_t^{\perp}$,
the first map is zero,
so $M_t^{-}\simeq M\oplus M_t^+[1]$.
Since $M_t^+[1]\in\mathcal U_t$ and $M_t^-\in\mathcal V_t$,
orthogonality gives $M_t^+=0$.
So $M\simeq M_t^{-}\in\mathcal{V}_t$.

\smallskip
\noindent\emph{Step 4:}
By \eqref{eq:calibration-confinement},
\[
\mathcal U_t\subseteq\cQ((t-\delta,\infty)),
\qquad
\mathcal V_t\subseteq\cQ((-\infty,t+\delta)).
\]
First fix $s<t$ with $t-s\leq\delta$ and take $0\ne F\in\mathcal U_s\cap\mathcal V_t$. 
Then
\[
F\in\cQ((s-\delta,t+\delta)).
\]
Since $t-s\leq\delta$,
applying the buffered HN construction in $\cQ((s-3\delta,t+5\delta))$,
we obtain a finite filtration of $F$ with factors $G_j\in\cR(\psi_j)$,
where $\psi_1>\cdots>\psi_\ell$.
The first and last HN factors give nonzero maps $G_1\to F$ and $F\to G_\ell$.
Since $F\in\mathcal V_t\cap\mathcal U_s$,
the Hom-vanishing \eqref{eq:calibration-hom} implies $\psi_1\leq t$ and $\psi_\ell>s$.
Consequently,
\[
s<\psi_j\leq t\qquad(1\leq j\leq\ell).
\]

Now let $0\ne M\in\DD$.
Choose real numbers $t_0<t_n$ satisfying
\[
t_0<\phi_{\cQ}^{-}(M)-\delta,
\qquad
t_n>\phi_{\cQ}^{+}(M)+\delta.
\]
By \eqref{eq:calibration-cut-bounds},
$M\in\mathcal U_{t_0}\cap\mathcal V_{t_n}$.
Choose a finite subdivision $t_0<t_1<\cdots<t_n$ with $t_{j+1}-t_j\leq\delta$.

Let $T_j\to M$ be the $\mathcal U_{t_j}$-truncation of $M$,
grouping these truncations together gives canonical maps
\[
0=T_n\longrightarrow T_{n-1} \longrightarrow\cdots\longrightarrow T_0=M.
\]
Set $C_j=\Cone(T_{j+1}\to T_j)$,
we have
\[
C_j\in\mathcal U_{t_j}\cap\mathcal V_{t_{j+1}}.
\]
Each nonzero $C_j$ therefore has a finite HN filtration with factors in $\cR(\psi)$ of phases in $(t_j,t_{j+1}]$.
Refining the successive truncation triangles by these filtrations gives a finite HN filtration of $M$ with strictly decreasing phases.

It remains to prove that $\cR(\psi)$ is extension-closed.
Consider a triangle $F_1\to N\to F_2$ with $F_1,F_2\in\cR(\psi)$.
By \eqref{eq:calibration-confinement},
$N\in\cQ((\psi-\delta,\psi+\delta))$.
The buffered HN construction in $\cQ((\psi-3\delta,\psi+5\delta))$ gives factors $G_j\in\cR(\psi_j)$ with $\psi_1>\cdots>\psi_\ell$.

If $\psi_1>\psi$,
then \eqref{eq:calibration-hom} gives $\Hom(G_1,F_1)=\Hom(G_1,F_2)=0$,
hence $\Hom(G_1,N)=0$,
contradicting the nonzero map $G_1\to N$.
Dually,
if $\psi_\ell<\psi$,
then $\Hom(F_1,G_\ell)=\Hom(F_2,G_\ell)=0$,
contradicting the nonzero quotient map $N\to G_\ell$.
Thus there is only one factor of phase $\psi$.
Hence $N\in\cR(\psi)$.

In conclusion,
we know $\cR$ is a slicing.
By its definition,
every nonzero $F\in\cR(\psi)$ satisfies $W(F)=|W(F)|e^{i\pi\psi}$.
Therefore $\cR$ is a slicing compatible with $W$.

\smallskip
\noindent\emph{Step 5:}
By \eqref{eq:calibration-confinement},
\[
\cR((t-\delta,t+\delta))
\subseteq
\cQ((t-2\delta,t+2\delta)).
\]
The category on the right is a thin quasi-abelian category of finite length.
Strict exact sequences on the left are also strict exact on the right,
so $\cR$ is locally finite and $(W,\cR)$ is a pre-stability condition.

Finally,
\eqref{eq:calibration-confinement} and \cite[Lemma~6.1]{Bri} give $d(\cQ,\cR)\leq\delta$.
\end{proof}

Suppose that the subcategories $U^{-1}\cP(\phi)$ form a slicing compatible with $W$. 
Since $U$ is exact and conservative, the induced slicing has coaisles $U^{-1}\cP(\leq c)$, which coincide with those of $\cQ$ by Proposition~\ref{prop:finite-transitions}.
Hence the two slicings coincide, and in particular $\cQ$ is compatible with $W$. 
Since Theorem~4.1 gives
$d(\cQ,\cR)\le\delta<1$ and both slicing are compatible with $W$, 
\cite[Lemma~6.4]{Bri}
implies that $\cR=\cQ$.
Thus our construction agrees with exact induction whenever the latter applies.

\section{The stability condition and full numerical support}
\label{sec:numerical}

\subsection{Construction and a mass estimate}
We now establish the support property of the pre-stability condition constructed in Section~\ref{sec:calibration} with respect to the full numerical Grothendieck group.

Fix \(\eps=\frac{1}{100}\) and $\delta=\frac{3}{2}\eps$. 
By Corollary~\ref{cor:induced-charge},
the hypotheses of Theorem~\ref{thm:calibration} are satisfied,
so we obtain a locally finite pre-stability condition $\tau=(W,\cR)$ on $\DD$,
where $W=Z\circ U$,
such that
\[
\cR(\psi)\subseteq
\cQ((\psi-\delta,\psi+\delta))
\qquad(\psi\in\R),
\]
and $d(\cQ,\cR)\leq\delta$.

The following lemma controls the mass of every $\tau$-semistable object.
It will be used to transfer the full numerical support property from $\sigma$ to $\tau$.

\begin{lemma}
\label{lem:mass estimate}
    Let $\tau=(W,\cR)$ be the pre-stability condition constructed in Section~\ref{sec:calibration}.
    Then for every nonzero $\tau$-semistable object $E\in\cR(\psi)$,
    one has
    \begin{equation}\label{eq: sigma phase of UE}
    UE\in\cP((\psi-\frac{5\eps}{2},\psi+\frac{3\eps}{2}))\subseteq\cP((\psi-3\eps,\psi+2\eps)).
    \end{equation}
    Consequently,
    the mass of $UE$ satisfies 
    \begin{equation}\label{eq: mass of UE}
     m_{\sigma}(UE)\leq \sec(3\pi\eps)|W(E)|.    
    \end{equation}
\end{lemma}
\begin{proof}
Let $0\neq E\in \cR(\psi)$.
By Proposition~\ref{prop:finite-transitions} and Theorem~\ref{thm:calibration},
one has 
$$
E\in\cQ((\psi-\delta,\psi+\delta)),\qquad U(\cQ(\phi))\subseteq \cP([\phi-\eps,\phi]),
$$
Since $U$ is exact,
applying to $\cQ$-HN filtration of $E$ gives $UE\in\cP((\psi-\delta-\eps,\psi+\delta))\subseteq\cP((\psi-\frac{5}{2}\eps,\psi+\frac{3}{2}\eps))$,
this proves \eqref{eq: sigma phase of UE}.
Let $F_1,\ldots,F_t$ be the $\sigma$-HN factors of $UE$ satisfying $F_j\in\cP(\phi_j)$.
By \eqref{eq: sigma phase of UE},
we know $|\phi_j-\psi|<3\eps$.
Since $W(E)=|W(E)|e^{i\pi\psi}$,
we obtain
$$
\begin{aligned}
|W(E)|=\operatorname{Re}\bigl(e^{-i\pi\psi}Z(UE)\bigr)
&=\sum_j |Z(F_j)|\cos\bigl(\pi(\phi_j-\psi)\bigr)\\
&\geq\cos(3\pi\eps)\sum_j|Z(F_j)|\\
&=\cos(3\pi\eps)m_\sigma(UE).  
\end{aligned}
$$
As $3\eps<\frac{1}{2}$,
$\cos(3\pi\eps)>0$,
and \eqref{eq: mass of UE} follows. 
\end{proof}

\subsection{Comparison of numerical Grothendieck groups}
To deduce the full numerical support property from Lemma~\ref{lem:mass estimate},
we need to compare the numerical Grothendieck groups of $\CC=\Db(\Coh(X))$ and $\DD=\Db(\Coh(X,\cA))$.
We will show that the forgetful functor $U$ induces an injection on numerical Grothendieck groups,
and it becomes an isomorphism after tensoring with $\Q$ or $\R$.
This allows us to pull the full-support estimate from $\sigma$ on $\Knum(\CC)$ back to $\tau$ on $\Knum(\DD)$.

Recall from Section~\ref{sec:cuts} that $\cA$ is an Azumaya algebra of degree $r$,
so $\rk_{\OO_X}\cA=r^2$. 
The free--forgetful adjunction is $$
L:\CC\rightleftarrows\DD:U,\qquad L(F)=\cA\otimes_{\OO_X} F,\qquad L\dashv U
$$
Both functors are exact and induce homomorphisms between Grothendieck groups of $\CC$ and $\DD$,
denoted by $L_*:K_0(\CC)\to K_0(\DD)$ and $U_*:K_0(\DD)\to K_0(\CC)$.
Since $X$ is smooth,
we identify $K_0(\CC)$ with the Grothendieck ring $K_0(X)$.

\begin{proposition}\label{lem:K0-comparison}
The induced homomorphisms \(L_*\) and \(U_*\) satisfy
\begin{equation}\label{eq: L_* and U_*}
L_*U_*=[\cA]\cdot \id_{K_0(\DD)},\qquad U_*L_*=[\cA]\cdot\id_{K_0(\CC)}    
\end{equation}
where the first equality uses the natural $K_0(\CC)$-module structure on $K_0(\DD)$.
The class $[\cA]$ is invertible in $K_0(\CC)_{\Q}$.
Consequently,
$U_*$ induces an isomorphism $U_{\Q}:K_0(\DD)_{\Q} \xrightarrow{\sim}K_0(\CC)_{\Q}.$ Moreover,
$U_*$ descends to an injective homomorphism
\[
U_{\mathrm{num}}:
\Knum(\DD)
\hookrightarrow \Knum(\CC),
\]
which becomes an isomorphism after tensoring with $\Q$ (and hence with $\R$).
\end{proposition}

\begin{proof}
We first prove \eqref{eq: L_* and U_*}.
For any $F\in\CC$,
the equality \(U_*L_*(F)=[\cA]\otimes F\) is immediate.
For the other equality,
let $M\in\DD$. 
The underlying $\OO_X$-complex $\cA\otimes_{\OO_X}UM$ carries two natural left $\cA$-module structures,
the first one is $b\cdot_1(a\otimes m)=ba\otimes m$ and the second is $b\cdot_2(a\otimes m)=a\otimes bm$.
The first is the free module structure on $LU(M)$,
while the second represents $[\cA]\cdot [M]$ under the natural $K_0(\CC)$-module structure on $K_0(\DD)$.

Now we show that these two structures are naturally isomorphic. 
Let $ \trd:\cA\longrightarrow\OO_X $ be the reduced trace.
It induces a nondegenerate pairing
\[
\cA\otimes_{\OO_X}\cA
\longrightarrow\OO_X,
\qquad
a\otimes b\longmapsto\trd(ab).
\]
Indeed,
on an \'etale splitting cover $V\to X$,
the pairing becomes the ordinary matrix-trace pairing.

Let
\[
\Omega=\sum_s u_s\otimes v_s
\in\Gamma(X,\cA\otimes_{\OO_X}\cA)
\]
be the canonical inverse tensor of the reduced-trace pairing.
On an \'etale splitting chart,
$\Omega=\sum_{i,j}E_{ij}\otimes E_{ji}.$ Direct computation gives
\begin{equation}
\Omega^2=1\otimes 1,\qquad \Omega(b\otimes 1)=(1\otimes b)\Omega
\end{equation}
for every local section $b$ of $\cA$.
The tensor $\Omega$ is canonical,
so these identities hold globally.
Using $\Omega$,
we define an $\OO_X$-linear morphism locally in an \'etale trivialization by 
$$
T_M:a\otimes m\longrightarrow\sum_s u_sa\otimes v_sm.
$$
The local expressions glue since $\Omega$ is canonical.
Then $T_M$ is an $\cA$-linear involution intertwining the two $\cA$-module structures,
in particular,
$T_M$ is an isomorphism.
This proves both identities in \eqref{eq: L_* and U_*}.

Next we show that $\left[\cA\right]$ is invertible in $K_0(X)_{\Q}$.
Write $d=\dim X$,
since $X$ is smooth,
the Chern character map induces an isomorphism of rings 
$$
\ch:K_0(X)_{\Q}\stackrel{\sim}{\longrightarrow}\bigoplus_{j=0}^d\operatorname{CH}^{j}(X)_{\Q}.
$$
As $X$ is connected,
the rank map is the degree-zero component of the Chern character.
Set $\xi=[\cA]-r^2[\OO_X]\in K_0(X)_{\Q}$,
then $\ch(\xi)\in \bigoplus_{j=1}^{d}\operatorname{CH}^j(X)_{\Q}$.
Since $\operatorname{CH}^{j}(X)_{\Q}=0$ for $j\geq d+1$,
then $[\cA]=r^2[\OO_X]+\xi$
is invertible,
with inverse $$
[\cA]^{-1}=\frac{1}{r^2}\sum_{j=0}^d\left(-\frac{\xi}{r^2} \right)^j.
$$
All inversions here take place in $K_0(\CC)_{\Q}$,
independent of the characteristic $k$.
Both $L_*$ and $U_*$ are compatible with $K_0(\CC)$-module structures, 
then \eqref{eq: L_* and U_*} implies $U_{\Q}:K_0(\DD)_{\Q}\to K_0(\CC)_{\Q}$ is an isomorphism with inverse $U_{\Q}^{-1}=L_{\Q}\circ([\cA]^{-1}\cdot-)$.

We now pass to numerical Grothendieck groups.
Let \(\chi_{\CC}\) and \(\chi_{\DD}\) denote the Euler pairings on $K_0(\CC)$ and $K_0(\DD)$.
Then derived adjunction $L\dashv U$ gives
\begin{equation}\label{eq:Euler-adjunction}
\chi_{\CC}(F,UM)=\chi_{\DD}(LF,M)
\end{equation}
for all $F\in \CC$ and $M\in\DD$.
The numerical Grothendieck groups are defined by $$
\Knum(\CC)=K_0(\CC)/\operatorname{rad}_R(\chi_{\CC}),\qquad \Knum(\DD)=K_0(\DD)/\operatorname{rad}_R(\chi_{\DD}).
$$
Direct computation shows $x\in\operatorname{rad}_R(\chi_{\DD})$ implies $U_*x\in\operatorname{rad}_R(\chi_{\CC})$,
so $U$ descends to a homomorphism $$
U_{\mathrm{num}}:\Knum(\DD)\longrightarrow \Knum(\CC).
$$
To prove injectivity,
suppose $x\in K_0(\DD)$ such that $U_*x\in\operatorname{rad}_R(\chi_{\CC})$,
then for an arbitrary $y\in K_0(\DD)$, 
there exists $f\in K_0(X)_{\Q}$ such that $y\otimes 1=L_{\Q}(f)$ since $L_{\Q}$ is surjective.
Extending the Euler pairing $\Q$-linearly,
we obtain $$
\chi_{\DD}(y,x)=\chi_{\DD}(L_{\Q}(f),x)=\chi_{\CC}(f,U_*x)=0.
$$
So $x\in \operatorname{rad}_R(\chi_{\DD})$,
this proves that $U_{\mathrm{num}}$ is injective.
Tensoring with $\Q$,
surjectivity of $U_{\Q}$ implies $U_{\mathrm{num},\Q}$ is surjective, 
together with injectivity of $U_{\mathrm{num}}$, 
this proves that $U_{\mathrm{num},\Q}$ is an isomorphism.

Finally,
$\Knum(\CC)$ is a free abelian group of finite rank,
as $X$ is smooth projective. 
The injectivity of $U_{\mathrm{num}}$ implies that $\Knum(\DD)$ is also free abelian of finite rank, 
tensoring $U_{\mathrm{num},\Q}$ with $\R$, 
we know $U_{\mathrm{num},\R}:\Knum(\DD)_{\R}\stackrel{\sim}{\longrightarrow}\Knum(\CC)_{\R}$.
This completes the proof.
\end{proof}

\begin{remark}
    The map $U_{\mathrm{num}}$ need not be an integral isomorphism,
    even for a split Azumaya algebra.
    
    Indeed,
    take $X=\operatorname{Spec}k$ and $\cA=\operatorname{Mat}_r(k)$ with $r>1$. 
    Let $S=k^r$ be the simple left $\cA$-module.
    Then we know $K_0(\DD)=\Z[S]$ and $K_0(k)=\Z[k]$,
    so forgetful functor satisfies $U_*[S]=r[k]$ and $L_*[k]=[\cA]=r[S]$,
    so $U_*L_*=L_*U_*=r^2$.
    This result also shows that the forgetful functor does not preserve the Euler pairing in general.
\end{remark}
We are now ready to complete the proof of Theorem~\ref{thm:main}.
\begin{proof}[Proof of Theorem~\ref{thm:main}]

First assume that $X$ is connected. 
The preceding construction gives the locally finite pre-stability condition \(\tau=(W,\cR)\) on $\DD$,
with $W=Z\circ U$.

Since the central charge $Z$ constructed in Section~\ref{sec:uniform} is numerical,
Proposition~\ref{lem:K0-comparison} proves that $\tau$ is also numerical,
since its central charge $W$ factors through $\Knum(\DD)$.

By \cite[Corollary~2.3]{Cheng},
the same stability condition \(\sigma=(Z,\cP)\) satisfies the support property with respect to \(\Knum(\CC)\).
Choose a Euclidean norm $\|\cdot\|_{\CC}$ on \(\Knum(\CC)_{\R}\) and a constant \(C_\sigma>0\) with \(\|[F]\|_{\CC}\le C_\sigma|Z(F)|\) for every nonzero $\sigma$-semistable object \(F\).

By Proposition~\ref{lem:K0-comparison},
\[
\|v\|_{\DD}
:=\|U_{\mathrm{num},\R}(v)\|_{\CC}
\]
defines a Euclidean norm on $\Knum(\DD)_{\R}$.
For a nonzero \(\tau\)-semistable object \(M\),
let \(F_1,\ldots, F_s\) be the $\sigma$-HN factors of \(UM\),
then
\begin{equation}
    \begin{aligned}
     \|[M]\|_{\DD} &=\|[UM]\|_{\CC}
\le\sum_j\|[F_j]\|_{\CC}\\
&\le C_\sigma\sum_j|Z(F_j)| \le C_\sigma\sec(3\pi\eps)|W(M)|,  
    \end{aligned}
\end{equation}
where the last inequality follows from Lemma~\ref{lem:mass estimate},
and $C_{\sigma}$ is independent of $M$.
Hence $\tau$ is a Bridgeland stability condition satisfying the full numerical support property.

Write $X=\bigsqcup_{i=1}^sX_i$ for the decomposition of $X$ into connected components.
Then $\Db(X,\alpha)=\bigoplus_{i=1}^s\Db(X_i,\alpha_i)$,
where $\alpha_i$ is the restriction of $\alpha$ to $X_i$.
Construct $\tau_i=(W_i,\cR_i)$ for each $i$ and define 
$$W(M_1,\ldots,M_s)=\sum_{i=1}^sW_i(M_i),\qquad \cR(\psi)=\bigoplus_{i=1}^s\cR_i(\psi)$$ it is easy to verify $(W,\cR)$ is the stability condition we want.
Finally, Morita equivalence sends this stability condition on $\Db(X,\cA)$ to $\Db(X,\alpha)$,
preserving the full numerical support property.
This completes the proof.
\end{proof}

\endgroup

\section{Properties of the constructed stability condition}
\label{sec:properties}

Throughout this section,
let $X$ be a smooth projective variety over an arbitrary field $k$,
and let $\cA$ be an Azumaya algebra representing $\alpha\in\operatorname{Br}(X)$.
We retain the notation
\[
\CC=\Db(\Coh X),\qquad
\DD=\Db(\Coh(X,\cA)),\qquad \sigma=(Z,\cP),\qquad \tau=(W,\cR)
\] from Sections~\ref{sec:cuts} to \ref{sec:numerical},
where $L\dashv U$ is the free-forgetful adjunction,
$W=Z\circ U$,
$\eps=\frac{1}{100}$ and $\delta=\frac{3\eps}{2}$.
We first assume that $X$ is connected,
so that $\cA$ has constant degree $r$.
The conclusions extend to disconnected $X$ by the same argument in Section~\ref{sec:numerical},
and hence to $\Db(X,\alpha)$ via Morita equivalence.

\subsection{Geometricity of $\tau$}\label{subsec:geometricity}
Let $K/k$ be a field extension.
Retain the notation $X_K$,
$\cA_K$ and $\DD_K$ from Section~\ref{subsec:prelim-stability} and put $\CC_K=\Db(\Coh(X_K))$.
By the scalar-extension theorem recalled in Section~\ref{subsec:prelim-stability},
the stability condition $\tau=(W,\cR)$ admits a canonical scalar extension $\tau_{K,1}=(W_{K,1},\cR_{K,1})$ to $\DD_K$.

On the other hand,
applying the construction of Section~\ref{sec:uniform}--\ref{sec:numerical} directly over $K$, 
using the base-changed projective embedding and the same parameter $M$ in Section~\ref{sec:uniform},
gives another stability condition $\tau_{K,2}=(W_{K,2},\cR_{K,2})$.

To describe the central charge of the canonical extension,
set $ N_{k}=\Knum(\DD)$,
$N_K=\Knum(\DD_K).
$ Recall the normalized Euler-test map from \cite[(12.2)]{BLMNPS},
\[
\eta^\vee_{K/k}:N_K
\longrightarrow N_k\otimes_{\Z}\Q,
\]
characterized by $\chi_k(T,\eta^\vee_{K/k}([E]))=\chi_K(T_K,E)$ for every $T\in\DD$ and $E\in\DD_K$.
Here $\chi_k$ and $\chi_K$ denote the Euler pairings on $\DD$ and $\DD_K$ respectively.

Since $W$ is numerical,
it admits a $\Q$-linear extension $W_{\Q}:N_k\otimes_{\Z}\Q\to\C$.
The central charge of the canonical scalar extension is therefore
\[
W_{K,1}=W_{\Q}\circ\eta^\vee_{K/k}.
\]

We now show that the canonical scalar extension agrees with the direct construction.
\begin{proposition}
\label{lem:properties-base-change}
$\tau_{K,1}=\tau_{K,2}$ as above.
\end{proposition}

\begin{proof}
Let $p:X_K\to X$ be the base-change projection,
which induces the restriction functor $\DD_{\mathrm{qc}}(X_K,\cA_K)\to \DD_{\mathrm{qc}}(X,\cA)$.
Let $\sigma_{K,1}=(Z_{K,1},\cP_{K,1})$ be the base change of $\sigma=(Z,\cP)$ by $p$ to $X_K$ and let $\sigma_{K,2}$ be the direct construction on $X_K$.
By \cite[Theorem~2.8 and Lemma~3.10]{LLLMPSZ},
pullback and pushforward commute with scalar extension,
hence $\sigma_{K,1}=\sigma_{K,2}$.

Moreover,
the two cutoff inclusions following from
\eqref{eq:fixed-error} extend to the Ind categories, written as 
$$
\cA\otimes \widehat{\cP}(>c)\subseteq\widehat{\cP}(>c-\eps),\qquad \cA\otimes\widehat{\cP}(\leq c)\subseteq \widehat{\cP}(\leq c+\eps).
$$
Together with the projection formula,
we know $U_KL_K(\cP_K(\phi))\subset \cP_K([\phi-\eps,\phi+\eps])$,
this is the base change version of equation \eqref{eq:fixed-error}.

By Lemma~\ref{lem:large-bands} and \cite[Theorem~2.8(1)]{LLLMPSZ},
we know 
\[
\begin{aligned}
\BB_{K,c}
&=\{E\in\DD_K:
      U_KE\in\cP_K(\le c)\}\\
&=\{E\in\DD_K:
      Up_*E\in\widehat{\cP}(\le c)\}\\
&=\{E\in\DD_K:
      p_*E\in\widehat{\BB}_c\},
\end{aligned}
\]
where $p_*U_K=Up_*$.

By Proposition~\ref{prop:coherence},
$\Ind(\cA_c)$ is an aisle.
Since $\cA_c=\widehat{\cA}_c\cap\DD$ and $\widehat{\cA}_c$ is generated by $L(\cP(>c))\subseteq \cA_c$, 
the two-sided inclusion gives $\widehat{\cA}_c=\Ind(\cA_c)$.
By Proposition~\ref{prop:finite-transitions},
we know $$
\widehat{\BB}_c=\widehat{\cQ}(\leq c),\qquad\widehat{\BB}_{K,c}=\widehat{\cQ}_K(\leq c).
$$
Consequently,
$\cQ_K(\leq c)=p_*^{-1}\widehat{\cQ}(\leq c)= \{E\in\DD_K:
p_*E\in\widehat{\cQ}(\le c)\}$.

The key step is the comparison of two slicings $\cR_{K,1}$ and $\cR_{K,2}$,
in fact,
we need to show $d(\cR_{K,1},\cR_{K,2})\leq 2\delta$.
Since $d(\cQ,\cR)\leq\delta$,
after taking Ind-extension,
we know
\begin{equation}
\widehat{\cR}(\leq c)
\subseteq\widehat{\cQ}(\leq c+\delta),
\qquad
\widehat{\cQ}(\leq c)
\subseteq\widehat{\cR}(\leq c+\delta).    
\end{equation}

By definition,
we know $\cR_{K,1}(\leq c)=\left\{E\in\DD_K:p_*E\in\widehat{\cR}(\leq c) \right\}$,
taking inverse images under $p_*$ and restricting to perfect objects gives
\begin{equation}
\cR_{K,1}(\leq c)
\subseteq\cQ_K(\leq c+\delta),\qquad \cQ_K(\leq c)
\subseteq\cR_{K,1}(\leq c+\delta).   
\end{equation}
Therefore $d(\cQ_K,\cR_{K,1})\leq\delta$,
on the other hand,
Theorem~\ref{thm:calibration} gives $d(\cQ_K,\cR_{K,2})\leq\delta$.
Consequently,
\[
d(\cR_{K,1},\cR_{K,2})
\leq2\delta<1.
\]

Finally,
$\tau_{K,1}$ and $\tau_{K,2}$ have the same central charge.
Indeed,
nondegeneracy of the rational numerical Euler pairing allows us to write
\[
Z(-)=\sum_\nu z_\nu\chi_{\CC}(T_\nu,-),
\qquad T_\nu\in K_0(\CC)_{\Q}.
\]
Since $L\dashv U$,
we know $W(-)=\sum_\nu z_\nu\chi_{\DD}(LT_\nu,-)$.
The coefficients $z_\nu$ are unchanged under scalar extension, 
and $(LT_\nu)_K=L_K((T_\nu)_K)$ by flat base change,
so
$$
\begin{aligned}
W_{K,1}(-)
&=\sum_\nu z_\nu
  \chi_{\DD_K}((LT_\nu)_K,-)\\
&=\sum_\nu z_\nu
  \chi_{\DD_K}(L_K(T_{\nu})_K,-)=Z_K(U_K(-))
=W_{K,2}(-).    
\end{aligned}
$$
Then \cite[Lemma~6.4]{Bri} gives $\tau_{K,1}=\tau_{K,2}$,
that finishes the proof.
\end{proof}

For a closed point $x\in X$,
let $F_x=\kappa(x)$ be its residue field and write
\[
\cA_x:=\cA\otimes_{\OO_X}\kappa(x)\simeq\operatorname{Mat}_{m_x}(D_x),\qquad d_x=\deg D_x,
\]
where $D_x$ is a central division algebra over $F_x$,
direct computation implies $r=m_xd_x$.
Write $S_x$ for its unique simple left module,
regarded as a coherent $\cA$-module on $X$ supported at $x$.
Then
\begin{equation}\label{eq:properties-point-modules}
U(S_x)\simeq F_x^{\oplus rd_x},\qquad
L(F_x)\simeq S_x^{\oplus m_x}.
\end{equation}
We call a condition on $\DD$ \emph{geometric} if,
after extension to an algebraic closure $\bar k$,
all simple twisted skyscrapers $S_p$ at geometric points $p\in X_{\overline{k}}$ are stable of the same phase.
This is the twisted analogue of \cite[Definition~2.11]{LLLMPSZ}.

\begin{proposition}\label{prop:properties-geometric}
The stability condition $\tau$ is geometric.
Moreover,
for every closed point $x\in X$,
the object $S_x$ is $\tau$-stable of phase $1$ over $k$.
\end{proposition}

\begin{proof}
We first check that $F_x$ is $\sigma$-stable of phase $1$.
Let $p:X_{\overline{k}}\to X$,
after pullback by $p^*$,
we know $p^*F_x\simeq \OO_{x\times_k\overline{k}}$,
where $x\times _k\overline{k}=\operatorname{Spec}(\kappa(x)\otimes_k\overline{k})$.
It has a finite filtration by geometric skyscraper sheaves of length $[F_x:k]$,
by \cite[Theorem~6.2]{LLLMPSZ}, 
each factor is semistable of phase 1,
lies in $\cP_{\overline{k}}(1)$,
so $p^*F_x\in\cP_{\overline{k}}(1)$.
By the construction of $\cP_{\overline{k}}$,
we know $F_x\in\cP(1)$.

To prove $F_x$ is stable,
for any short exact sequence $0\to A\to F_x\to B\to 0$ in $\cP(1)$,
then \cite[Theorem~2.8(3)]{LLLMPSZ} implies that after pull back by $p^*$,
\begin{equation}\label{eq:p^* pullback}
0\longrightarrow p^*A\longrightarrow p^*F_x\longrightarrow p^*B\to 0    
\end{equation}
is a short exact sequence in $\cP_{\overline{k}}(1)$.
Since $p^*F_x$ has finite length in $\cP_{\bar k}(1)$ and all its simple composition factors are geometric skyscraper sheaves, 
the same holds for $p^*A$ and $p^*B$.
Thus both are extensions of zero-dimensional coherent sheaves and hence belong to $\Coh(X_{\bar k})$.
As $p$ is faithfully flat,
we obtain $A,B\in\Coh(X)$.
The simplicity of $F_x\in\operatorname{Coh(X)}$ concludes the proof.
The definition of central charge \eqref{eq: formula of central charge} gives $Z(F_x)=-[F_x:k]$.

Now we turn to the twisted category. 
For every $a<1$,
$F_x\in\cP(1)\subseteq\cP(>a)$,
we know $L(F_x)\in\bigcap_{a<1}\cA_a$.
On the other hand,
$UL(F_x)\simeq F_x^{\oplus r^2}\in \cP(\leq 1)$. 
So by Lemma~\ref{lem:large-bands},
we know $L(F_x)\in\BB_1\subset\bigcap_{b>1}\BB_b$.
Thus,
by the definition of $\cQ(1)$,
$L(F_x)\in\cQ(1)$,
means $S_x\in\cQ(1)$ since $L(F_x)\simeq S_x^{\oplus m_x}$ and $\cQ(1)$ is closed under direct summands. 
Moreover,
$W(S_x)=-rd_x[F_x:k]$ by direct computation. 

Set $h=2\delta$,
$I=(1-h,1+h)$ and $\mathcal{E}=\cQ(I)$,
so $S_x$ is enveloped by $\mathcal{E}$.
Consider the strict exact sequence
\begin{equation}\label{eq:properties-point-subobject}
0\longrightarrow F\longrightarrow S_x\longrightarrow G
\longrightarrow0
\end{equation}
in $\cQ(I)$.
By \eqref{eq:slicing-phase-range},
all three underlying objects lie in
\[
\cP((1-h-\eps,1+h))\subset\cP((1/2,3/2]).
\]
The latter is an abelian heart,
so $U$ induces a short exact sequence $0\to UF\to F_x^{\oplus rd_x}\to UG\to 0$ in $\HH$,
thus $\phi_{\cP}^+(UF)\leq \phi_{\cP}^+(U(S_x))=1$.
If $F\neq0$,
then $UF\neq0$,
and the argument of $W(F)=Z(UF)$ is at most $\pi$.
This proves $W$-semistability in the thin enveloping category,
hence $S_x\in\cR(1)$ by the definition of $(W,\cR)$.

Finally we show that $S_x$ is $\tau$-stable.
For the same short exact sequence \eqref{eq:properties-point-subobject},
we assume $\phi_W(F)=1$ and $F\neq 0$.
Since $\phi_W(F)=1$,
we know $W(F)=Z(UF)\in\R_{<0}$.
Moreover,
by the preceding phase bound $UF\in \cP((1-h-\eps,1])$.
Thus the imaginary parts of the central charges of its HN-factors are nonnegative,
consequently,
all these factors are of phase $1$,
so $UF\in\cP(1)$.
Then $UG\in\cP([1,1+h))$,
and since $W(G)=W(S_x)-W(F)$,
so $W(G)\in\R$,
that means $G=0$ or $UG\in\cP(1)$.
Since $F_x$ is simple in $\cP(1)$ and $UF$ is a subobject of $U(S_x)\simeq F_x^{\oplus rd_x}$ in $\cP(1)$,
then $UF\simeq F_x^{\oplus e}$ for some $1\leq e\leq rd_x$.

Since $U$ is $t$-exact for the standard $t$-structures and $UF\simeq F_x^{\oplus e}$,
we know $H^j(F)= 0$ for $j\neq 0$ and then $F$ is a finite-dimensional left $\cA_x$-module.
By the classification of $\cA_x$-modules,
we know $F\simeq S_x^{\oplus q}$,
then $UF\simeq F_x^{\oplus qrd_x}$,
so $q=1$ since $1\leq e\leq rd_x$,
then $F\to S_x$ is an isomorphism,
which forces $G=0$.

Finally,
every short exact sequence in $\cR(1)$ has its terms in $\cQ((1-\delta,1+\delta))\subset\cQ(I)$ by
\eqref{eq:calibration-confinement}, and is strict in $\cQ(I)$, so $S_x$ is simple in $\cR(1)$.
Applying the same argument over $\bar k$ and using Proposition~\ref{lem:properties-base-change},
we know $\tau$ is geometric.
\end{proof}

\begin{remark}
The free fiber $L(F_x)$ is stable exactly when $r/d_x=1$,
which is equivalent to $\cA_x$ being a division algebra,
and otherwise is strictly semistable.

But in general,
$\tau$-stability of $S_x$ over $k$ does not imply that  $(S_x)_{\overline{k}}$-stable over $\overline{k}$.
For example,
let $k=\R$,
$X=\operatorname{Spec}(\R)$,
and $\cA=\mathbb H$,
the Hamilton quaternion algebra. 
Its unique simple left module is $S=\mathbb H$.

After scalar extension to $\C$,
one has
\[
\mathbb H\otimes_{\R}\C
\simeq\operatorname{Mat}_2(\C).
\]
The base-changed module $S_{\C}$ is the regular left $\operatorname{Mat}_2(\C)$-module,
hence is the direct sum of two copies of its unique simple module.
It is therefore strictly semistable.
\end{remark}

\subsection{The mass--Hom bound}
A stability condition $\rho$ on a proper $k$-linear dg-category $\DD$ satisfies a \textit{mass--Hom bound} if, 
for every $T\in\DD$,
there exists $C_T>0$ such that
\[
\dim_k\Hom_{\DD}(T,E) \leq C_Tm_\rho(E)
\]
for all $E\in\DD$.
Since mass is invariant under shifts,
the same constant $C_T$ applies to $E[j]$ for every $j\in\Z$.

\begin{lemma}\label{lem:properties-all-object-mass}
For every $E\in\DD$,
\begin{equation}\label{eq:properties-all-object-mass}
m_\sigma(UE)\leq\sec(3\pi\eps) m_\tau(E).
\end{equation}
\end{lemma}

\begin{proof}
Let $E_1,\ldots,E_t$ be the $\tau$-HN factors of $E$.
Applying $U$ to the HN triangles of $E$,
then Lemma~\ref{lem:mass estimate} gives
\[
m_\sigma(UE)
\leq\sum_{j=1}^t m_\sigma(UE_j)
\leq \sec(3\pi\eps)\sum_{j=1}^t |W(E_j)|
=\sec(3\pi\eps) m_\tau(E).
\]
\end{proof}

\begin{proposition}\label{prop:properties-mass-Hom}
For every $T\in\DD$ there is a constant $C_T>0$ such that
\begin{equation}\label{eq:properties-mass-Hom}
\dim_k\Hom_{\DD}(T,E[j])\leq C_T m_\tau(E)
\qquad \text{for every }E\in\DD,\ j\in\Z.
\end{equation}
In particular,
$\tau$ satisfies the mass--Hom bound.
\end{proposition}

\begin{proof}
Choose a compact generator $G$ of $\Dqc(X)$,
since $X$ is smooth projective,
we know $D_{\mathrm{qc}}(X)^{c}=\Perf(X)=\CC$,
so $G\in\CC$.
Then we know $\CC=\operatorname{thick}(G)$ and $\DD=\operatorname{thick}(LG)$.
By \cite[Theorem~7.5]{LLLMPSZ},
the stability condition $\sigma$ has a mass--Hom bound.
Thus there is a constant $B_G>0$ satisfying
\[
\dim_k\Hom_{\CC}(G,F)\leq B_Gm_\sigma(F)
\qquad \text{for every }F\in\CC.
\]
For any $E\in\DD$ and $j\in\Z$,
adjunction and Lemma~\ref{lem:properties-all-object-mass} yield
\begin{equation}\label{eq: mass-Hom bound of LG}
\begin{aligned}
\dim_k\Hom_{\DD}(LG,E[j])
&=\dim_k\Hom_{\CC}(G,UE[j])\\
&\leq B_Gm_\sigma(UE[j])
\leq B_G\sec(3\pi\eps)\, m_\tau(E).
\end{aligned}    
\end{equation}
Here mass is invariant under shifts,
so the same constant works for every integer $j$.

Let $$
\mathcal{T}:=\left\{ T\in \DD:\exists\, C_T>0 \text{ satisfying \eqref{eq:properties-mass-Hom} for all $E$, $j$}\right\}.$$
It is easy to show $\mathcal{T}$ is closed under shifts,
finite direct sums,
direct summands and cones,
so $\mathcal{T}$ is a thick triangulated subcategory.
And \eqref{eq: mass-Hom bound of LG} implies $LG\in \mathcal{T}$,
so $\mathcal{T}=\DD$ since $\DD=\operatorname{thick}(LG)$. 
That completes the proof.
\end{proof}

The mass--Hom bound also persists throughout the connected component containing $\tau$.
Put $N=\Knum(\DD)$,
and write $\Stab_N(\DD)$ for the space of numerical Bridgeland stability conditions whose central charges factor through $N$ and satisfy the support property with respect to $N$.

\begin{corollary}\label{cor:properties-component-mass-hom}
Let $\Sigma$ be the connected component of $\Stab_N(\DD)$ containing $\tau$.
Then every $\rho\in\Sigma$ satisfies a mass--Hom bound.
\end{corollary}

\begin{proof}
By Lemma~\ref{prelim:topology},
for every $\rho\in\Sigma$ there exists $C_\rho>0$ such that
\[
m_\tau(E)\leq C_\rho m_\rho(E)
\qquad \text{for every $E\in\DD$.}
\]
For each $T\in\DD$,
Proposition~
\ref{prop:properties-mass-Hom} therefore gives
\[
\dim_k\Hom_{\DD}(T,E[j]) \leq C_T C_\rho\,m_\rho(E)
\qquad \text{for every $E\in\DD,\ j\in\Z$.}
\]
Hence $\rho$ satisfies the mass--Hom bound.
\end{proof}

\bibliographystyle{amsalpha}
\bibliography{reference}

\end{document}